\documentclass[11pt]{article}
\usepackage[T1]{fontenc}
\usepackage{lmodern,amsmath,amssymb,amsthm,mathtools}
\usepackage[margin=1in]{geometry}
\usepackage{microtype,tikz-cd}
\usepackage{array}
\usepackage[numbers]{natbib}
\usepackage[colorlinks=true,linkcolor=blue,citecolor=blue,urlcolor=blue]{hyperref}
\hypersetup{pdftitle={Maximality of the Dosen--Petric Boolean Category (v27)},pdfauthor={Florian Lengyel}}

\newtheorem{theorem}{Theorem}[section]
\newtheorem{lemma}[theorem]{Lemma}
\newtheorem{proposition}[theorem]{Proposition}
\newtheorem{corollary}[theorem]{Corollary}
\newcommand{\B}{\mathsf B}
\newcommand{\Occ}{\operatorname{Occ}}
\newcommand{\Rel}{\mathbf{Rel}_{\mathrm{fin}}}
\newcommand{\FormulaSet}{\operatorname{Form}}
\title{Maximality of the Do\v{s}en--Petri\'c Boolean Category}
\author{Florian Lengyel\thanks{Email: \href{mailto:flengyel@gradcenter.cuny.edu}{\texttt{flengyel@gradcenter.cuny.edu}}}}
\date{November 17, 2026\quad(v27)}
\begin{document}
\maketitle

\begin{abstract}
We prove that the zero-identity Boolean category of Do\v{s}en and
Petri\'c, freely generated by a nonempty set of propositional variables,
is maximal. Every congruence preserved by the operations on morphisms
and by uniform formula substitution either is equality or identifies
every pair of parallel morphisms. This disproves their conjecture that
the category is not maximal. We give two proofs that use Boolean
Coherence and uniform formula substitution to derive equality of the
product projections from an identification of distinct parallel
morphisms. We also give a Lean formalization of these proofs and their
dependencies. It includes Boolean Coherence and the correspondence
between substitution on tagged occurrence sets and substitution on
finite ordinals.
\end{abstract}

\medskip
\noindent\textit{Mathematics Subject Classification (2020).}
Primary 03G30; Secondary 03F07, 18A32.

\smallskip
\noindent\textit{Keywords.} Boolean categories, maximality, categorial proof theory,
coherence, substitution congruences.

\section{Introduction}

We prove that the freely generated zero-identity Boolean category of
Do\v{s}en and Petri\'c over a nonempty set of propositional variables
is maximal. Adding any equation between parallel morphisms that does
not follow from its defining equations, and closing under the operations on
morphisms and uniform formula substitution, collapses the category
to a preorder. This disproves Do\v{s}en and Petri\'c's conjecture that the category
is not maximal \cite[\S14.2, p.~322]{DP}. Before adjoining the equation, the
two projections $p\wedge p\to p$ are distinct for every
propositional variable $p$.

The proof reduces an added equation to equality of the product
projections, as in the maximality theorem for cartesian categories
\cite[\S5]{DPCartesian}. Boolean Coherence supplies a pair of
variable occurrences belonging to one morphism's relation and not the
other. Precomposition and postcomposition isolate this pair and
give a literal instance of equation~\eqref{eq:context} in
Section~\ref{sec:isolation} \cite[\S9.7, pp.~215--216]{DP}.
Uniform substitution and the double-negation isomorphism then give
this equation for every formula.
The isomorphisms of Lemma~\ref{lem:structural-isos}, whose types are
listed in \cite[\S12.5, p.~294]{DP}, give equality of the product
projections. Sections~\ref{sec:isolation}
and~\ref{sec:collapse} present this proof; Section~\ref{sec:formalization}
gives a second proof using mix. The proof of
\cite[Lemma~7]{DPIsomorphic} also supplies the two isomorphisms.

Other notions of Boolean category need not make conjunction and
disjunction categorical products and coproducts. Free categories constructed from classical proof nets
are studied in \cite[Theorem~3.5.4]{LSFree}, equations on monoid
and comonoid structures in $*$-autonomous categories in
\cite{Strassburger}, and inequalities interpreting cut reduction
in \cite{FP}. Section~\ref{sec:category} compares the negation
equation of \cite{FP} with the zero-identity equation used here.

We use the occurrence coordinates of \cite[\S9]{Notation}:
occurrences form tagged disjoint unions, and substitution replaces
each occurrence by a block. Section~\ref{sec:occurrences} defines
these coordinates; Appendix~\ref{app:substitution} proves their
correspondence with the finite-ordinal substitution formula of
Do\v{s}en and Petri\'c. Section~\ref{sec:formalization} and
Appendix~\ref{app:lean-statements} describe the Lean formalization,
including a guide to its declarations and proof steps.
Appendix~\ref{app:explicit-inverses} constructs the two inverses
used in the first maximality proof, and
Appendix~\ref{app:cut-errata} corrects two typographical errors
in the cut-elimination argument used in Do\v{s}en and Petri\'c's
proof of Boolean Coherence \cite{DP}.
References to \cite{DP} use the printed pagination of its
September 2007 revised version. Lean declaration names in footnotes
omit the namespace prefix \texttt{DosenPetric.}

\section{The Boolean proof category}\label{sec:category}

Fix a nonempty set $P$ of propositional variables. Write $\FormulaSet(P)$
for the formulas built from $P$, the constants $\top,\bot$, the binary
connectives $\wedge,\vee$, and negation $\neg$. Let $\B(P)$ be the
freely generated zero-identity Boolean category defined in
\cite[\S\S14.2--14.3]{DP}. Its objects are these formulas, and its
morphisms are the arrow terms of that presentation modulo its equations.
We write $\B$ when $P$ is fixed. We recall the structure and equations
used in the proof.

The connectives $\wedge$ and $\vee$ give binary products and coproducts;
$\top$ is terminal and $\bot$ is initial. For $A_1=A$ and $A_2=B$,
denote the projections and injections by
\[
 \pi_i^{A,B}:A\wedge B\longrightarrow A_i,
 \qquad
 \iota_i^{A,B}:A_i\longrightarrow A\vee B
 \qquad(i=1,2).
\]
For $f:X\to A$ and $g:X\to B$, write
$\langle f,g\rangle:X\to A\wedge B$ for their pairing.
For $f:A\to X$ and $g:B\to X$, write
$[f,g]:A\vee B\to X$ for their copairing. Write
$!_A:A\to\top$ and $\kappa_A:\bot\to A$ for the terminal and initial
morphisms, and $1_A:A\to A$ for the identity. Composition $g\circ f$
applies $f$ first. For $f:A\to A'$ and $g:C\to C'$, the product and
coproduct bifunctors give morphisms
\[
 f\wedge g:A\wedge C\longrightarrow A'\wedge C',
 \qquad f\vee g:A\vee C\longrightarrow A'\vee C'.
\]

The dissociativity morphism is
\[
 d_{A,B,C}:A\wedge(B\vee C)\longrightarrow(A\wedge B)\vee C.
\]
Its first-injection equation is
\begin{equation}
 d_{A,B,C}\circ(1_A\wedge\iota_1^{B,C})
   =\iota_1^{A\wedge B,C}.
 \label{eq:first-injection}
\end{equation}
This is the coproduct-injection equation of
\cite[\S11.1, p.~234]{DP}.

Negation is represented by a functor
$\neg:\B^{\mathrm{op}}\to\B$, sending $f:A\to C$ to
$\neg f:\neg C\to\neg A$. Its double-negation and De Morgan
isomorphisms are natural. In particular, write
$n_A:\neg\neg A\to A$ for the double-negation isomorphism.
The presentation also has morphisms
$\eta_A:\top\to\neg A\vee A$ and
$\epsilon_A:A\wedge\neg A\to\bot$.
The \emph{zero identity} is the endomorphism $0_A:A\to A$ defined by
\begin{equation}
 0_A=[\kappa_A,1_A]\circ(\epsilon_A\vee1_A)\circ d_{A,\neg A,A}
     \circ(1_A\wedge\eta_A)\circ\langle1_A,!_A\rangle.
 \label{eq:zero-identity}
\end{equation}
The defining zero-identity equation asserts that
\begin{equation}
 f\circ0_A=0_C\circ g
 \qquad(f,g:A\to C).
 \label{eq:zero-law}
\end{equation}
Equation~\eqref{eq:zero-identity} is $(0\eta\epsilon)$ of
\cite[\S14.2, p.~317]{DP}. Its product and coproduct unit isomorphisms
are expressed here in pairing and copairing notation.
Equation~\eqref{eq:zero-law} is the zero-identity law $(0I)$;
see \cite[\S14.3, p.~322]{DP}.

Appendix~\ref{app:structural-diagrams} displays
equations~\eqref{eq:first-injection}--\eqref{eq:zero-law}
as commutative diagrams.

Whenever $\B(A,C)$ is nonempty, define $0_{A,C}=f\circ0_A$ for any
$f:A\to C$. Equation~\eqref{eq:zero-law} makes this independent of
$f$: both $f\circ0_A$ and $g\circ0_A$ equal $0_C\circ g$.
Taking $f=g$ gives $f\circ0_A=0_C\circ f$. Consequently zero morphisms
absorb composition whenever the composites are defined:
$h\circ0_{A,C}=0_{A,D}$ for $h:C\to D$, and
$0_{A,C}\circ k=0_{X,C}$ for $k:X\to A$.

F\"uhrmann and Pym interpret cut reduction by inequalities between
morphisms. An inequality belongs to a net theory if and only if it
holds in every interpretation of that theory in their classical
categories \cite[Theorems~6.4--6.5]{FP}.
Their negation triangle makes the composite in
\eqref{eq:zero-identity}, with the corresponding monoidal unit
isomorphisms, equal to $1_A$ \cite[(11)]{FP}.
Here that composite is $0_A$.

Two morphisms are \emph{parallel} when they have the same source and
target. A \emph{proof congruence} on $\B$ is an equivalence relation
on each hom-set preserved by composition and every operation on morphisms
of the presentation, including the derived operations of pairing,
copairing, and negation.
A \emph{substitution congruence} is a proof congruence also preserved
by every uniform substitution $\sigma:P\to\FormulaSet(P)$.
Here substitution acts on the entire arrow term, including its source
and target; this is rule~$(su)$ of \cite[\S2.7, p.~51]{DP}.

The \emph{total congruence} identifies every parallel pair. A category
is \emph{thin}, or a \emph{preorder}, when each hom-set has at most one
element. For a proof congruence $\theta$, the quotient $\B/\theta$
has the same objects as $\B$ and morphisms given by $\theta$-equivalence
classes. It is a preorder exactly when $\theta$ is total.

\begin{theorem}\label{thm:main}
Let $P$ be a nonempty set of propositional variables, and let
$\B(P)$ be the freely generated zero-identity Boolean category of
Do\v{s}en and Petri\'c. For every substitution congruence $\theta$
on $\B(P)$ that identifies two distinct parallel morphisms, the
quotient $\B(P)/\theta$ is a preorder.
\end{theorem}

\emph{Maximality} means that this collapse occurs for every added
equation between parallel arrow terms that does not follow from the
defining equations, with closure under the operations on morphisms and
uniform substitution. Equivalently, the only substitution
congruences are equality and the total congruence.
Sections~\ref{sec:isolation} and~\ref{sec:collapse} give the proof.

\section{Occurrences and substitution}\label{sec:occurrences}

For a formula $A$, write $\Occ(A)$ for the finite set of positions
occupied by propositional variables in $A$. Repeated occurrences of the
same variable occupy distinct positions. We use disjoint unions to
distinguish these positions.

For finite sets $X,Y$, let
$X\sqcup Y=(\{1\}\times X)\cup(\{2\}\times Y)$ be their tagged
disjoint union. The occurrence sets are defined recursively by
\[
 \begin{aligned}
 \Occ(p)&=\{*\},&
 \Occ(\top)&=\Occ(\bot)=\varnothing,\\
 \Occ(\neg A)&=\Occ(A),&
 \Occ(A\mathbin\star C)&=\Occ(A)\sqcup\Occ(C)
       \quad(\star\in\{\wedge,\vee\}).
 \end{aligned}
\]
The label map $\lambda_A:\Occ(A)\to P$ assigns to each occurrence
its underlying variable. Distinct positions may have the same label.
Numbering the occurrences from left to right gives a bijection
$\Occ(A)\to\{1,\ldots,|\Occ(A)|\}$; the latter set is empty when
$|\Occ(A)|=0$. Appendix~\ref{app:ordinals} gives the substitution
formula in these finite-ordinal coordinates.

Let $\Rel$ be the category of finite sets and binary relations.
A relation $R:X\to Y$ is a subset of $X\times Y$; identities are
diagonal relations, and composition is relational composition.
Transporting the occurrence interpretation of \cite{DP} to the
tagged coordinates gives a functor
$G:\B\to\Rel$ with $GA=\Occ(A)$. For $f:A\to C$, its image is
the relation $Gf\subseteq\Occ(A)\times\Occ(C)$.
For example,
\[
 G\pi_1^{A,C}=\{((1,x),x):x\in\Occ(A)\},\qquad
 G\iota_1^{A,C}=\{(x,(1,x)):x\in\Occ(A)\}.
\]
The second projection and injection use tag~$2$. Pairing takes the
union of the two relations after tagging their targets; copairing
takes the union after tagging their sources. The relation $Gd_{A,B,C}$
is the graph of the rebracketing bijection
\[
 \Occ(A)\sqcup\bigl(\Occ(B)\sqcup\Occ(C)\bigr)
 \longrightarrow
 \bigl(\Occ(A)\sqcup\Occ(B)\bigr)\sqcup\Occ(C),
\]
which preserves each occurrence of $A,B,C$
\cite[\S2.9, p.~60]{DP}. Zero morphisms have empty relations.

For an occurrence $x\in\Occ(A)$, define
$\operatorname{sgn}_A(x)=(-1)^k\in\{+1,-1\}$, where $k$ is the
number of negations above $x$ in $A$. Every morphism $f:A\to C$
preserves both the variable and its sign:
\[
 (x,y)\in Gf\quad\Longrightarrow\quad
 \lambda_A(x)=\lambda_C(y),\qquad
 \operatorname{sgn}_A(x)=\operatorname{sgn}_C(y).
\]
Induction on arrow terms proves this assertion. Identities,
associativity and symmetry maps, projections and injections,
diagonals and codiagonals, mix, dissociativity and unit maps
preserve the variable and its sign. The double-negation and
De Morgan maps preserve the signs, and the zero,
$\eta$ and $\epsilon$ maps have empty relations. Composition and
the two binary operations preserve the assertion.\footnote{\raggedright Lean:
\href{https://github.com/flengyel/Dosen-Petric-Maximality/blob/dosen-petric-v25-palomar-v3/DosenPetric/PresentationSignedLabels.lean\#L73}{\texttt{Presentation.\allowbreak{}occurrence\_signedLabel}}.}

\begin{theorem}[Boolean Coherence, Do\v{s}en--Petri\'c]
\label{thm:coherence}
The functor $G:\B\to\Rel$ is faithful: parallel morphisms $f,g$ of
$\B$ are equal if and only if $Gf=Gg$.
\end{theorem}

This is \cite[\S14.2, p.~321]{DP}. In particular, for $p\in P$, the
two projections $p\wedge p\to p$ have different occurrence
relations, so $\B$ is nonthin.

Joyal's collapse result says that a cartesian closed category with
an initial object $\bot$ is a preorder if
$A\cong(A\Rightarrow\bot)\Rightarrow\bot$ for every object $A$
\cite[\S14.3, p.~324]{DP}.
The category $\B$ is not cartesian closed. Faithfulness gives at most
one morphism $\top\to E$ for every $E$, since $\Occ(\top)$ is empty,
whereas $1_p\ne0_p$ for $p\in P$. An exponential $p\Rightarrow p$
would give a bijection $\B(p,p)\cong\B(\top,p\Rightarrow p)$,
a contradiction.

The maximality theorem for cartesian closed categories
\cite[\S11]{DPCCC} uses a typed B\"ohm theorem.
The absence of $p\Rightarrow p$ prevents its application here.

A \emph{literal} is $p$ or $\neg p$ for $p\in P$.
A formula is in \emph{negation normal form} when negation occurs only
immediately before a propositional variable. The double-negation and
De Morgan isomorphisms, including the isomorphisms for the truth
constants, give an isomorphism $a_A:A\to A^\circ$ to a formula in
negation normal form. The relation $Ga_A$ is the graph of a bijection
$\Occ(A)\to\Occ(A^\circ)$
\cite[\S14.1, pp.~313--315]{DP}.

For morphisms between formulas in negation normal form, every pair in
the occurrence relation consists of occurrences of the same literal.
In the presentation on negation normal forms in
\cite[\S14.2, pp.~318--321]{DP}, objects are formulas built from
literals, conjunction, disjunction and the truth constants.
Arrow terms are generated by identities, associativity and symmetry
maps, diagonals and codiagonals, projections and injections, mix,
dissociativity, unit maps, zero endomorphisms, and $\eta_p,\epsilon_p$.
The zero endomorphisms and $\eta_p,\epsilon_p$ have empty relations;
each remaining generator joins equal literals.
Composition, $\wedge$ and $\vee$ preserve this property.
The equivalence with $\B$ preserves the occurrence interpretation,
so the assertion follows.

Negation takes the converse relation:
\[
 (y,x)\in G(\neg f)\quad\Longleftrightarrow\quad(x,y)\in Gf
 \qquad(f:A\to C).
\]
On the normal-form presentation, negation interchanges the two
binary operations and reverses each generator relation and the
order of composition. Induction on terms therefore gives the
converse relation. The occurrence bijections of normalization
transport this equation to the original positions, using
$\Occ(\neg A)=\Occ(A)$ and $\Occ(\neg C)=\Occ(C)$.\footnote{\raggedright Lean:
\href{https://github.com/flengyel/Dosen-Petric-Maximality/blob/dosen-petric-v25-palomar-v3/DosenPetric/ReducedDualRelations.lean\#L20}{\texttt{Reduced.Term.\allowbreak{}dual\_relation}} and
\href{https://github.com/flengyel/Dosen-Petric-Maximality/blob/dosen-petric-v25-palomar-v3/DosenPetric/PresentationNegationOccurrence.lean\#L25}{\texttt{Presentation.\allowbreak{}occurrence\_negate}}.}

Substitution replaces each occurrence of a variable by a separate
block of occurrences of its substituted formula.
Appendix~\ref{app:profiles} describes these blocks and their label maps
by commutative diagrams, following \cite[\S9]{Notation}.
Let $Q$ be a set of propositional variables and
$\sigma:P\to\FormulaSet(Q)$ a substitution. If $A[\sigma]$ denotes the
result of uniform substitution, there is a bijection
\begin{equation}
 \Occ(A[\sigma])\cong D_{A,\sigma},\qquad
 D_{A,\sigma}=\coprod_{x\in\Occ(A)}
                  \Occ\bigl(\sigma(\lambda_A(x))\bigr).
 \label{eq:occ-substitution}
\end{equation}
The set on the right consists of pairs $(x,t)$, where $x$ is an
occurrence of $A$ and $t$ is an occurrence within the formula inserted there.
Induction on $A$ proves the bijection: the variable case has one
block, the constant cases have none, binary connectives form tagged
unions, and negation preserves the occurrence set. In particular,
two occurrences of the same variable give two separately tagged
copies of its substituted formula.

Put $E_\sigma=\coprod_{p\in P}\Occ(\sigma(p))$.
Define $\operatorname{pr}_D:D_{A,\sigma}\to\Occ(A)$ and
$\operatorname{pr}_E:E_\sigma\to P$ by
$\operatorname{pr}_D(x,t)=x$ and $\operatorname{pr}_E(p,t)=p$.
The function
$\widehat\lambda_{A,\sigma}:D_{A,\sigma}\to E_\sigma$ sends
$(x,t)$ to $(\lambda_A(x),t)$. Thus
\begin{equation}
\begin{tikzcd}[column sep=large,row sep=large]
 D_{A,\sigma}\arrow[r,"\widehat\lambda_{A,\sigma}"]
   \arrow[d,"\operatorname{pr}_D"']
   &E_\sigma\arrow[d,"\operatorname{pr}_E"]\\
 \Occ(A)\arrow[r,"\lambda_A"']&P.
\end{tikzcd}
\label{eq:block-square}
\end{equation}
Both composites send $(x,t)$ to $\lambda_A(x)$.
For $(x,t),(x',t)\in D_{A,\sigma}$, their images under
$\widehat\lambda_{A,\sigma}$ agree exactly when
$\lambda_A(x)=\lambda_A(x')$, even if $x\ne x'$.

The same presentation defines $\B(Q)$ when $Q$ is empty.
Substitution on arrow terms induces a functor
\[
 F_\sigma:\B(P)\longrightarrow\B(Q),\qquad
 A\longmapsto A[\sigma],\quad f\longmapsto f[\sigma].
\]
Here $f:A\to C$ is sent to
$f[\sigma]:A[\sigma]\to C[\sigma]$.
Rule~$(su)$ makes this functor well defined on arrow terms modulo the
defining equations. It preserves identities, composition, and the
specified operations on morphisms. For each morphism, the blocks in
\eqref{eq:occ-substitution} describe the occurrence sets of its
substituted source and target.

Write $\chi_{A,\sigma}:\Occ(A[\sigma])\to D_{A,\sigma}$ for
the bijection in \eqref{eq:occ-substitution}. For $f:A\to C$,
define the relation
$R_{f,\sigma}:D_{A,\sigma}\to D_{C,\sigma}$ by
\[
 \begin{aligned}
 ((x,t),(y,u))\in R_{f,\sigma}
 \quad\Longleftrightarrow\quad&
 (x,y)\in Gf\ \text{ and}\\
 &(\lambda_A(x),t)=(\lambda_C(y),u)\text{ in }E_\sigma.
 \end{aligned}
\]
Thus an edge of $Gf$ is replaced by the identity relation within
the inserted block. The substituted relation satisfies
\[
 (x',y')\in G(f[\sigma])
 \quad\Longleftrightarrow\quad
 (\chi_{A,\sigma}(x'),\chi_{C,\sigma}(y'))\in R_{f,\sigma}.
\]
Regard the two bijections as their graphs. The following square
then commutes in $\Rel$:
\[
\begin{tikzcd}[column sep=large,row sep=large]
 \Occ(A[\sigma])\arrow[r,"{G(f[\sigma])}"]
    \arrow[d,"{\chi_{A,\sigma}}"']
   &\Occ(C[\sigma])\arrow[d,"{\chi_{C,\sigma}}"]\\
 D_{A,\sigma}\arrow[r,"{R_{f,\sigma}}"']&D_{C,\sigma}.
\end{tikzcd}
\]
Induction on the arrow term proves the relation formula.
Identities, associativity and symmetry maps, projections and
injections, diagonals and codiagonals, mix, dissociativity and
unit maps act on the inserted blocks by the same tagged relations.
The double-negation and De Morgan maps give their occurrence
bijections on these blocks; zero, $\eta$ and $\epsilon$ give
empty relations. The two binary
operations retain the tags. For composition, an intermediate old
occurrence supplies an intermediate block, and preservation of
variable labels identifies its inner occurrence on both sides.
The formula passes to quotient morphisms because both interpretations
respect the defining equations.\footnote{\raggedright Lean:
\href{https://github.com/flengyel/Dosen-Petric-Maximality/blob/dosen-petric-v25-palomar-v3/DosenPetric/PresentationSubstitutionRelation.lean\#L61}{\texttt{Presentation.\allowbreak{}occurrence\_substitute\_expanded}};
the composition step is
\href{https://github.com/flengyel/Dosen-Petric-Maximality/blob/dosen-petric-v25-palomar-v3/DosenPetric/SubstitutionFiberRelations.lean\#L48}{\texttt{Relations.Fiber.\allowbreak{}lift\_composition}}.}
Repeated old occurrences retain their different first coordinates;
an empty inserted block supplies no positions.

Let $R$ be a set of variables and let
$\tau:Q\to\FormulaSet(R)$ be another substitution. Define
$\tau\star\sigma:P\to\FormulaSet(R)$ by
$(\tau\star\sigma)(p)=\sigma(p)[\tau]$. For formulas $A,C$ over
$P$ and $f:A\to C$, substitution satisfies
\[
 \begin{aligned}
 A[p\mapsto p]&=A,&
 (A[\sigma])[\tau]&=A[\tau\star\sigma],\\
 f[p\mapsto p]&=f,&
 (f[\sigma])[\tau]&=f[\tau\star\sigma].
 \end{aligned}
\]
The morphism equalities use the displayed formula equalities to
identify their sources and targets. Pointwise equal substitutions
give equal substituted formulas and morphisms with these same
identifications. Formula induction proves the object equations;
term induction and passage to the quotient prove the morphism
equations.\footnote{\raggedright Lean:
\href{https://github.com/flengyel/Dosen-Petric-Maximality/blob/dosen-petric-v25-palomar-v3/DosenPetric/FormulaSubstitutionLaws.lean\#L11}{\texttt{Formula.\allowbreak{}subst\_identity}},
\href{https://github.com/flengyel/Dosen-Petric-Maximality/blob/dosen-petric-v25-palomar-v3/DosenPetric/FormulaSubstitutionLaws.lean\#L20}{\texttt{Formula.\allowbreak{}subst\_composition}}, and
\href{https://github.com/flengyel/Dosen-Petric-Maximality/blob/dosen-petric-v25-palomar-v3/DosenPetric/PresentationSubstitutionActionQuotient.lean\#L39}{\texttt{Presentation.\allowbreak{}substitute\_identity}},
\href{https://github.com/flengyel/Dosen-Petric-Maximality/blob/dosen-petric-v25-palomar-v3/DosenPetric/PresentationSubstitutionActionQuotient.lean\#L44}{\texttt{Presentation.\allowbreak{}substitute\_composition}},
\href{https://github.com/flengyel/Dosen-Petric-Maximality/blob/dosen-petric-v25-palomar-v3/DosenPetric/PresentationSubstitutionActionQuotient.lean\#L57}{\texttt{Presentation.\allowbreak{}substitute\_ext}}.}

\section{Isolating an occurrence}\label{sec:isolation}

Put $j_A=\iota_1^{A,\top}:A\to A\vee\top$. We shall derive the equation
\begin{equation}
 j_A\circ\pi_1^{A,\bot}
   =0_{A\wedge\bot,A\vee\top}.
 \tag{$K_A$}\label{eq:context}
\end{equation}
The right-hand side is also
$\iota_2^{A,\top}\circ!_\bot\circ\pi_2^{A,\bot}$: its occurrence
relation is empty, so Boolean Coherence identifies it with the zero
morphism. Thus \eqref{eq:context} is the equation used in
\cite[\S9.7, pp.~215--216]{DP}.

For negation-free formulas $A,B$ such that $A\to B$ is a
tautology, \cite[Lemma~6]{DPIsomorphic} constructs a morphism whose
occurrence relation consists of the edge between any prescribed
occurrences of the same variable in $A$ and $B$. The following
lemma isolates an edge of either of two given morphisms; normalization
allows the formulas to contain negation.

Selecting a pair on which relations differ also occurs in the
maximality proofs for presentations of split preorders, split
equivalence relations and binary relations \cite[\S16]{DPSplit}.

\begin{lemma}[Isolation of an occurrence]\label{lem:isolation}
For distinct parallel morphisms $f,g:A\to C$ in $\B$, there are a
literal $\ell$ and morphisms
\[
 i:\ell\wedge\bot\longrightarrow A,
 \qquad r:C\longrightarrow\ell\vee\top
\]
such that, after possibly interchanging $f$ and $g$,
\[
 r\circ f\circ i=j_\ell\circ\pi_1^{\ell,\bot},
 \qquad
 r\circ g\circ i=0_{\ell\wedge\bot,\ell\vee\top}.
\]
\end{lemma}

\begin{proof}
First replace $f,g$ by $a_C\circ f\circ a_A^{-1}$ and
$a_C\circ g\circ a_A^{-1}$. These remain distinct, so it suffices to
construct the morphisms for the negation normal forms $A^\circ,C^\circ$.
If
$i^\circ:\ell\wedge\bot\to A^\circ$ and
$r^\circ:C^\circ\to\ell\vee\top$ satisfy the equations,
then $i=a_A^{-1}\circ i^\circ:\ell\wedge\bot\to A$ and
$r=r^\circ\circ a_C:C\to\ell\vee\top$ satisfy them for the original
formulas.
We may therefore assume that $A,C$ are in negation normal form.

Faithfulness of $G$ gives a pair belonging to exactly one of $Gf,Gg$.
Interchange $f,g$ if necessary, and choose
$(x,y)\in Gf\setminus Gg$. The occurrences $x$ and $y$ have the same
literal $\ell$. Set $D=\ell\wedge\bot$ and $E=\ell\vee\top$.
For every formula $V$, the morphisms
\[
 b_V=\kappa_V\circ\pi_2^{\ell,\bot}:D\to V,
 \qquad
 t_V=\iota_2^{\ell,\top}\circ!_V:V\to E
\]
have empty occurrence relations.

Construct $i:D\to A$ by induction on $A$, following the subformula
that contains $x$ at each binary connective. At the selected literal, use $\pi_1^{\ell,\bot}$.
At a conjunction, pair the recursively constructed morphism with $b_V$
for the unselected factor, placing it in the order specified by the
tag. At a disjunction, compose the recursively constructed morphism
with the injection of the selected summand. The displayed rules for
$G$ show inductively that $Gi$ consists of one pair, whose first component
is the sole occurrence of $D$ and whose second component is $x$.

Construct $r:C\to E$ by induction on $C$, following the subformula
that contains $y$.
At the selected literal, use $\iota_1^{\ell,\top}$.
At a conjunction, first project to the selected factor and then apply
the recursively constructed morphism. At a disjunction, copair that
morphism with $t_V$ on the unselected summand. Thus $Gr$ consists of one pair,
whose first component is $y$ and whose second component is the sole
occurrence of $E$.

Relational composition now makes $G(r\circ f\circ i)$ the singleton
relation from $\Occ(D)$ to $\Occ(E)$, whereas
$G(r\circ g\circ i)$ is empty. Boolean Coherence gives the asserted
equations.
\end{proof}

The equations $(K_A)$ pass to retracts. Let $\theta$ be a proof
congruence on $\B$, and suppose that $i:A\to B$ and $r:B\to A$
satisfy $r\circ i=1_A$. If $(K_B)$ holds modulo $\theta$, then
$(K_A)$ holds modulo $\theta$.
Precomposing $(K_B)$ with $i\wedge1_\bot$ and postcomposing it
with $r\vee1_\top$ gives $(K_A)$ modulo the same proof
congruence. Indeed, the left side becomes
\[
 (r\vee1_\top)\circ j_B\circ\pi_1^{B,\bot}\circ(i\wedge1_\bot)
 =j_A\circ(r\circ i)\circ\pi_1^{A,\bot}
 =j_A\circ\pi_1^{A,\bot},
\]
and zero absorption gives the zero on the right.\footnote{\raggedright Lean:
\href{https://github.com/flengyel/Dosen-Petric-Maximality/blob/dosen-petric-v25-palomar-v3/DosenPetric/Substitution.lean\#L74}{\texttt{ContextualCollapse.\allowbreak{}context\_retract}}.}

\begin{proposition}\label{prop:all-contexts}
If a substitution congruence $\theta$ on $\B$ identifies distinct
parallel morphisms, then \eqref{eq:context} holds modulo $\theta$ for
every formula $A$.
\end{proposition}

\begin{proof}
Lemma~\ref{lem:isolation} and compatibility with composition give
\eqref{eq:context} for a literal $\ell$.
If $\ell=p$, let $\sigma$ substitute $A$ for $p$ and fix the other variables.
This yields \eqref{eq:context} for $A$: under
\eqref{eq:occ-substitution}, the literal occurrence in each endpoint
becomes the block $\Occ(A)$, and the operations on morphisms are preserved
by $F_\sigma$; see also Appendix~\ref{app:shared}.

If $\ell=\neg p$, substitute $\neg A$ for $p$, obtaining the equation
for $\neg\neg A$. Precompose with
$n_A^{-1}\wedge1_\bot$ and postcompose with $n_A\vee1_\top$,
where $n_A:\neg\neg A\to A$.
Put
$k_A=j_A\circ\pi_1^{A,\bot}:A\wedge\bot\to A\vee\top$.
The product and coproduct equations give the commutative square
\[
\begin{tikzcd}[column sep=large,row sep=large]
 A\wedge\bot\arrow[r,"k_A"]\arrow[d,"n_A^{-1}\wedge1_\bot"']
   &A\vee\top\\
 \neg\neg A\wedge\bot\arrow[r,"k_{\neg\neg A}"']
   &\neg\neg A\vee\top\arrow[u,"n_A\vee1_\top"'].
\end{tikzcd}
\]
Precomposition and postcomposition of the zero morphism give a zero
morphism by absorption. The equation for $\neg\neg A$ therefore gives
\eqref{eq:context} for $A$.
\end{proof}

\section{Equality of the product projections}\label{sec:collapse}

The next lemma identifies two isomorphisms of the types listed in
\cite[\S12.5, p.~294]{DP}.

\begin{lemma}\label{lem:structural-isos}
For every formula $A$, the morphisms
\[
 \begin{aligned}
 q_A&=[j_A\circ\pi_1^{A,\bot},\iota_2^{A,\top}]
      :(A\wedge\bot)\vee\top\longrightarrow A\vee\top,\\
 e_A&=1_A\wedge j_A
      :A\wedge A\longrightarrow A\wedge(A\vee\top)
 \end{aligned}
\]
are isomorphisms in $\B$.
\end{lemma}

\begin{proof}
First take $A=p\in P$. By \cite[\S12.5, p.~294]{DP}, the source
and target of $q_p$ are isomorphic. Each has one occurrence of $p$,
so the relation of any such isomorphism is the unique bijection
between these occurrence sets. The morphism $q_p$ has the same relation.
Boolean Coherence therefore identifies $q_p$ with that isomorphism.

The same cited list gives an isomorphism
$b:p\wedge p\to p\wedge(p\vee\top)$: take $n=m=1$ in
$p^n\wedge(p^m\vee\top)\cong p^{n+m}$ and choose $p^2=p\wedge p$.
The relation $Gb$ is a bijection between two-element occurrence sets.
After precomposing $b$ with the product symmetry if necessary, it
agrees with the order-preserving bijection $Ge_p$. Coherence then
identifies $e_p$ with an isomorphism. These two identifications also
follow from the occurrence-bijection construction in the proof of
\cite[Lemma~7]{DPIsomorphic}.

For arbitrary $A$, substitute $A$ for $p$ and fix the other variables.
The functor $F_\sigma$ sends $q_p,e_p$ to $q_A,e_A$ and preserves
their inverse equations.
\end{proof}

Appendix~\ref{app:explicit-inverses} constructs inverses for
these morphisms and verifies all four inverse equations.\footnote{\raggedright Lean:
\href{https://github.com/flengyel/Dosen-Petric-Maximality/blob/dosen-petric-v25-palomar-v3/DosenPetric/PresentationShortQ.lean\#L17}{\texttt{Presentation.ShortProof.\allowbreak{}qInverseTerm}} and
\href{https://github.com/flengyel/Dosen-Petric-Maximality/blob/dosen-petric-v25-palomar-v3/DosenPetric/PresentationShortE.lean\#L30}{\texttt{Presentation.ShortProof.\allowbreak{}eInverseTerm}}.}

\begin{lemma}\label{lem:collapse}
Let $\theta$ be a proof congruence on $\B$ such that
\eqref{eq:context} holds modulo $\theta$ for every formula $A$.
Then $\B/\theta$ is a preorder.
\end{lemma}

\begin{proof}
Write $\equiv$ for equality modulo $\theta$.
The quotient retains products and coproducts: their pairing,
copairing, and defining equations pass to equivalence classes.
Fix $A$, and put
$E=A\vee\top$, $S=(A\wedge\bot)\vee\top$, and
$s=\iota_2^{A,\top}:\top\to E$.
The morphisms $0_{A\wedge\bot,E}$ and $s\circ!_{A\wedge\bot}$ have
empty occurrence relations and hence are equal in $\B$. Thus
\eqref{eq:context} and the coproduct equations give
\[
 q_A\equiv[s\circ!_{A\wedge\bot},s]=s\circ!_S.
\]
The inverse equations for $q_A$ hold in the quotient, so
\[
 1_E=q_A\circ q_A^{-1}
    \equiv s\circ!_S\circ q_A^{-1}=s\circ!_E.
\]
Together with $!_E\circ s=1_\top$, this makes $E$ terminal in the
quotient.

Consequently the projection
$r=\pi_1^{A,E}:A\wedge E\to A$ is invertible, with inverse
$\langle1_A,s\circ!_A\rangle$.
Lemma~\ref{lem:structural-isos} makes $e_A$ invertible too, and the
product equations give
$r\circ e_A=\pi_1^{A,A}$.
Hence $\pi_1^{A,A}$ is invertible. Its right inverse is the diagonal
$\Delta_A=\langle1_A,1_A\rangle:A\to A\wedge A$, which must
therefore be its inverse. It follows that
\[
 \pi_2^{A,A}
   \equiv\pi_2^{A,A}\circ\Delta_A\circ\pi_1^{A,A}
   =\pi_1^{A,A}.
\]
For parallel morphisms $f,g:X\to A$, pairing gives
$f=\pi_1^{A,A}\circ\langle f,g\rangle
 \equiv\pi_2^{A,A}\circ\langle f,g\rangle=g$.
Thus the quotient is a preorder.
\end{proof}

\begin{proof}[Proof of Theorem~\ref{thm:main}]
Suppose that a substitution congruence $\theta$ on $\B(P)$ identifies
two distinct parallel morphisms. Proposition~\ref{prop:all-contexts}
gives \eqref{eq:context} modulo $\theta$ for every formula $A$.
Lemma~\ref{lem:collapse} then makes $\B(P)/\theta$ a preorder.
\end{proof}

We also obtain maximality in the sense of
\cite[\S9.3, p.~194]{DP}. An \emph{interpretation} in a zero-identity
Boolean category $\mathcal C$ assigns an object of $\mathcal C$ to
each propositional variable. For an assignment
$v:P\to\operatorname{Ob}(\mathcal C)$, the induced functor
$I_v:\B(P)\to\mathcal C$ preserves the specified products, coproducts,
units, dissociativity morphisms, negation maps, $\eta$ and $\epsilon$.

\begin{corollary}\label{cor:models}
Let $P$ be a nonempty set of propositional variables, and let
$f,g:A\to C$ be distinct morphisms of $\B(P)$.
If a zero-identity Boolean category $\mathcal C$ satisfies $f=g$
under every interpretation of the propositional variables, then
$\mathcal C$ is a preorder.
\end{corollary}

\begin{proof}
Declare two parallel morphisms of $\B(P)$ equivalent when all their
interpretations in $\mathcal C$ agree. Since every $I_v$ preserves
the operations on morphisms, this relation is a proof congruence. For a
substitution $\sigma:P\to\FormulaSet(P)$, the composite
$I_v\circ F_\sigma$ is the interpretation induced by
$p\mapsto I_v(\sigma(p))$. The congruence is therefore preserved
by substitution. It contains $f=g$, so
Theorem~\ref{thm:main} makes the quotient a preorder. Hence the
congruence is total. In particular,
$\pi_1^{p,p}=\pi_2^{p,p}$ holds under every interpretation, for any
$p\in P$. Interpret $p$ as an arbitrary object $X$ of $\mathcal C$.
Its two projections $X\wedge X\to X$ are equal. For any morphisms $a,b:Y\to X$, composing
$\langle a,b\rangle:Y\to X\wedge X$ with the equal projections gives $a=b$.
\end{proof}

For zero-identity distributive lattice categories,
\cite[\S12.5, p.~292]{DP} proves maximality without the truth
constants and relative maximality when both constants are present.
For lattice and dicartesian categories without an assumption of
distribution, the maximality and relative-maximality results are
proved in \cite[\S\S8--10]{DPBicartesian}.

\section{A second proof and the Lean formalization}\label{sec:formalization}

The inverse terms in Appendix~\ref{app:explicit-inverses} and the
coherence argument in Appendix~\ref{app:coherence-proof} are taken
from the Lean development. Appendix~\ref{app:reading-lean}
explains how the declarations and proof steps express the mathematical
statements used here. The footnotes identify the source declarations.

In the Lean~4 development \cite{Formalization}, formulas and arrow
terms are inductive types, with arrow terms indexed by their source
and target. The relation \texttt{Derives} is generated by the
presentation's equations and congruence rules. Its quotient is the
category $\B(P)$. Substitution, normalization and the occurrence
interpretation are defined on terms and proved to respect
\texttt{Derives}.

The theorem \texttt{relative\_maximality} states that a
substitution congruence identifying distinct parallel morphisms is total.
It assumes the operations on morphisms, a faithful occurrence interpretation,
isomorphisms to negation normal forms, preservation of literals
along occurrence relations, uniform substitution and double
negation, with their equations.
Its proof that the quotient is a preorder uses the morphism
\[
 h_A=[1_A,1_A]\circ m_{A,A}:A\wedge A\longrightarrow A,
\]
where the mix morphism $m_{A,C}:A\wedge C\to A\vee C$ is given, for any
formula $B$, by
\[
 m_{A,C}=(\pi_1^{A,B}\vee1_C)\circ d_{A,B,C}
              \circ(1_A\wedge\iota_2^{B,C}).
\]
This is equation $(dm)$ of \cite[\S11.1, p.~234]{DP}.
The right-hand side is independent of $B$. The symmetry of mix gives
$h_A\circ\langle\pi_2^{A,A},\pi_1^{A,A}\rangle=h_A$.
Write $\equiv$ for equality modulo a proof congruence containing \eqref{eq:context}
for every formula. For $B=\bot$, put
\[
 v_A=d_{A,\bot,A}\circ(1_A\wedge\iota_2^{\bot,A})
 :A\wedge A\longrightarrow(A\wedge\bot)\vee A.
\]
The coproduct equations and \eqref{eq:context} give
\[
 \begin{aligned}
 j_A\circ h_A
 &= [j_A\circ\pi_1^{A,\bot},j_A]\circ v_A\\
 &\equiv[0_{A\wedge\bot,A\vee\top},j_A]\circ v_A
  =j_A\circ\pi_2^{A,A}.
 \end{aligned}
\]
The last equality holds in $\B$ by Boolean Coherence: both occurrence
relations send the second copy of $\Occ(A)$ identically to the
$A$-summand of $\Occ(A\vee\top)$ and have no pairs from the first copy.
Precomposing with the symmetry of $A\wedge A$ and using the
preceding identity for $h_A$ gives
$j_A\circ\pi_1^{A,A}\equiv j_A\circ\pi_2^{A,A}$.
Precompose this equation with $\langle1_A,0_A\rangle$. Since the
composite $j_A\circ0_A$ has empty occurrence relation, coherence gives
\[
 j_A\equiv\iota_2^{A,\top}\circ!_A.
\]
Put $s_A=\iota_2^{A,\top}:\top\to A\vee\top$. The coproduct
equations give
\[
 1_{A\vee\top}=[j_A,s_A]
   \equiv[s_A\circ!_A,s_A]=s_A\circ!_{A\vee\top}.
\]
Since $!_{A\vee\top}\circ s_A=1_\top$, the object $A\vee\top$ is
terminal in the quotient. Composing with the coproduct symmetry gives
$\iota_2^{\top,A}\equiv\iota_1^{\top,A}\circ!_A$.
For $B=\top$, put
$w_A=1_A\wedge!_A:A\wedge A\to A\wedge\top$. Substituting the
preceding equation for $\iota_2^{\top,A}$ in the expression for mix gives
\[
 \begin{aligned}
 h_A&\equiv[\pi_1^{A,\top},1_A]\circ d_{A,\top,A}
             \circ(1_A\wedge\iota_1^{\top,A})\circ w_A\\
    &=[\pi_1^{A,\top},1_A]\circ\iota_1^{A\wedge\top,A}\circ w_A\\
    &=\pi_1^{A,\top}\circ w_A=\pi_1^{A,A}.
 \end{aligned}
\]
The second line uses \eqref{eq:first-injection}; the third uses the
coproduct and product equations. The symmetry of $h_A$ yields $\pi_1^{A,A}\equiv\pi_2^{A,A}$,
and pairing equates any two parallel morphisms.

For parallel morphisms $f,g:A\to C$, define
\[
 f\cup g=h_C\circ(f\wedge g)\circ\Delta_A:A\longrightarrow C,
 \qquad \Delta_A=\langle1_A,1_A\rangle.
\]
This is the union operation of \cite[p.~220]{DP}, written with
the product and coproduct operations of $\B$. Its occurrence
relation is
\[
 G(f\cup g)=Gf\cup Gg.
\]
Indeed, the diagonal supplies both source copies, $f\wedge g$
acts on the corresponding copies, and $h_C$ joins the two target
copies. Faithfulness gives, for parallel $f,g,k:A\to C$,
\[
 (f\cup g)\cup k=f\cup(g\cup k),\qquad
 f\cup g=g\cup f,\qquad f\cup f=f.
\]
For $a:X\to A$ and $b:C\to Y$, relational composition and
faithfulness give
\[
 \begin{aligned}
 (f\cup g)\circ a&=(f\circ a)\cup(g\circ a),\\
 b\circ(f\cup g)&=(b\circ f)\cup(b\circ g).
 \end{aligned}
\]
On an inhabited hom-set $\B(A,C)$, its zero satisfies
$f\cup0_{A,C}=f$; see also \cite[pp.~290--291]{DP}.
Consequently $\B(A,C)$ is a join semilattice with least element
$0_{A,C}$, where
\[
 f\leq g\quad\Longleftrightarrow\quad f\cup g=g
 \quad\Longleftrightarrow\quad Gf\subseteq Gg.
\]
Finite unions are independent of the order and repeated
summands. On each nonempty hom-set, the empty union is its zero morphism.\footnote{\raggedright The Lean laws are
\href{https://github.com/flengyel/Dosen-Petric-Maximality/blob/dosen-petric-v25-palomar-v3/DosenPetric/ReducedUnion.lean\#L131}{\texttt{Reduced.\allowbreak{}union\_assoc}},
\href{https://github.com/flengyel/Dosen-Petric-Maximality/blob/dosen-petric-v25-palomar-v3/DosenPetric/ReducedUnion.lean\#L90}{\texttt{Reduced.\allowbreak{}union\_comm}},
\href{https://github.com/flengyel/Dosen-Petric-Maximality/blob/dosen-petric-v25-palomar-v3/DosenPetric/ReducedUnion.lean\#L97}{\texttt{Reduced.\allowbreak{}union\_idem}},
\href{https://github.com/flengyel/Dosen-Petric-Maximality/blob/dosen-petric-v25-palomar-v3/DosenPetric/ReducedUnion.lean\#L80}{\texttt{Reduced.\allowbreak{}union\_pre}},
\href{https://github.com/flengyel/Dosen-Petric-Maximality/blob/dosen-petric-v25-palomar-v3/DosenPetric/ReducedUnion.lean\#L85}{\texttt{Reduced.\allowbreak{}union\_post}}, and
\href{https://github.com/flengyel/Dosen-Petric-Maximality/blob/dosen-petric-v25-palomar-v3/DosenPetric/ReducedFiniteUnion.lean\#L57}{\texttt{Reduced.finiteUnion.\allowbreak{}eq\_of\_members}}.
These Lean laws concern the presentation on negation normal forms
described in Section~\ref{sec:occurrences}. The mutually inverse maps
\href{https://github.com/flengyel/Dosen-Petric-Maximality/blob/dosen-petric-v25-palomar-v3/DosenPetric/ReducedQuotientReduction.lean\#L21}{\texttt{Reduced.\allowbreak{}reduceHom}} and
\href{https://github.com/flengyel/Dosen-Petric-Maximality/blob/dosen-petric-v25-palomar-v3/DosenPetric/ReducedRecovery.lean\#L26}{\texttt{Reduced.\allowbreak{}recoverHom}} identify its hom-set from $A^\circ$ to $C^\circ$ with
$\B(P)(A,C)$, transferring these laws to $\B(P)$.}

Cut elimination and Boolean Coherence are proved for the presentation
defining $\B(P)$. Normalization, preservation of literals and
uniform substitution on terms and their source and target formulas
satisfy the remaining hypotheses of \texttt{relative\_maximality}.
Its application gives \texttt{Presentation.maximality}.

The theorem \texttt{Presentation.short\_maximality} formalizes the proof
in Section~\ref{sec:collapse}, using the inverse terms for $q_A$ and $e_A$.
The model corollary is proved using interpretations preserving the
specified operations and their composition with substitution. The formalization
also verifies the occurrence maps and finite-ordinal constructions of
Appendix~\ref{app:substitution}, the diagram identities of
Appendix~\ref{app:structural-diagrams}, and the cut reductions and rank
comparisons discussed in Appendix~\ref{app:cut-errata}.

Five declarations are registered in the Palomar Registry
\cite{palomar-2026-10-30-000014-v3}:
the two maximality proofs, the congruence dichotomy, nonzero variable
identities, and the absence of morphisms from $\top$ to $\bot$.
Appendix~\ref{app:lean-statements} gives their statements in mathematical
notation and relates the arrow terms to morphisms of $\B(P)$.
It also describes the construction of the coherence proof.

\appendix
\section{Substitution in occurrence coordinates}\label{app:substitution}

Do\v{s}en and Petri\'c's generalized composition of definable
connectives \cite[\S2.5, pp.~45--46]{DP} substitutes formula
patterns and computes the map labelling the resulting occurrences.
We describe this map by the tagged coordinates of
\cite[\S\S9.2--9.3]{Notation}. Each occurrence after substitution
is indexed by an occurrence of the outer pattern and an occurrence
of the pattern inserted there.

\subsection{Occurrence maps and their fibres}\label{app:profiles}

For an integer $n\geq0$, put $\underline n=\{1,\ldots,n\}$, with
$\underline0=\varnothing$. A \emph{formula pattern} on a finite set
of placeholders is a formula in those placeholders. Its
\emph{occurrence map} records the index of the placeholder at each position,
after enumerating the placeholders and reading positions from left to right.

Let $U=\{p_1,\ldots,p_k\}$ consist of distinct placeholders and let
$M\in\FormulaSet(U)$ have $m$ occurrences. Its occurrence map is
\[
 o_M:\underline m\longrightarrow\underline k,
 \qquad o_M(r)=i\quad\text{when position $r$ contains $p_i$}.
\]
For each $i\in\underline k$, choose a set
$Q_i=\{q_{i,1},\ldots,q_{i,\ell_i}\}$ of distinct local placeholders and a
pattern $N_i\in\FormulaSet(Q_i)$ with $n_i$ occurrences. The sets
$Q_i$ are pairwise disjoint. Write
$o_i:\underline{n_i}\to\underline{\ell_i}$ for the occurrence map of $N_i$.
Empty placeholder sets and empty occurrence sets are allowed; a
pattern over an empty placeholder set has no variable occurrences.

Define the sets of occurrence coordinates and placeholder indices by
\[
 \begin{aligned}
 D&=\{(r,t):r\in\underline m,\ t\in\underline{n_{o_M(r)}}\},\\
 C&=\{(i,s):i\in\underline k,\ s\in\underline{\ell_i}\}.
 \end{aligned}
\]
Thus $(r,t)$ specifies an occurrence in the copy of $N_{o_M(r)}$
inserted at position $r$ of $M$. The pair $(i,s)$ specifies a local
placeholder of $N_i$. Define the injections
\[
 a_r:\underline{n_{o_M(r)}}\longrightarrow D,\quad a_r(t)=(r,t),
 \qquad
 b_i:\underline{\ell_i}\longrightarrow C,\quad b_i(s)=(i,s).
\]
There is a unique function $\gamma:D\to C$ for which, for every
$r\in\underline m$, the following square commutes:
\begin{equation}\label{eq:app-fibre}
\begin{tikzcd}[column sep=large,row sep=large]
 \underline{n_{o_M(r)}}\arrow[r,"a_r"]\arrow[d,"o_{o_M(r)}"']
   &D\arrow[d,"\gamma"]\\
 \underline{\ell_{o_M(r)}}\arrow[r,"b_{o_M(r)}"']&C.
\end{tikzcd}
\end{equation}
Every element of $D$ is $a_r(t)$ for a unique pair $(r,t)$.
Thus these squares determine $\gamma$ uniquely by
\[
 \gamma(r,t)=\bigl(o_M(r),o_{o_M(r)}(t)\bigr).
\]

Define the projections $p_D:D\to\underline m$ and
$p_C:C\to\underline k$ by $p_D(r,t)=r$ and $p_C(i,s)=i$.
Then $p_C\circ\gamma=o_M\circ p_D$, so the following square commutes:
\begin{equation}\label{eq:app-profile}
\begin{tikzcd}[column sep=large,row sep=large]
 D\arrow[r,"\gamma"]\arrow[d,"p_D"']&C\arrow[d,"p_C"]\\
 \underline m\arrow[r,"o_M"']&\underline k.
\end{tikzcd}
\end{equation}
Both composites send $(r,t)$ to $o_M(r)$.
If $r\ne r'$ but $o_M(r)=o_M(r')=i$, the copies of
$\underline{n_i}$ indexed by $r$ and $r'$ are disjoint subsets
of $D$.

To identify this construction with pattern substitution,
replace each $q_{i,s}$ in $N_i$ by $(i,s)$ and denote the resulting
pattern over $C$ by $\widetilde N_i$. Let
$\Phi_0:U\to\FormulaSet(C)$ be given by $\Phi_0(p_i)=\widetilde N_i$.
Every occurrence of $M[\Phi_0]$ has unique coordinates $(r,t)\in D$.
This gives a bijection $\chi:\Occ(M[\Phi_0])\to D$, and its relation
to the label map is the commutative square
\begin{equation}\label{eq:app-occurrence-labels}
\begin{tikzcd}[column sep=large,row sep=large]
 \Occ(M[\Phi_0])\arrow[r,"\chi","\cong"']
    \arrow[d,"\lambda_{M[\Phi_0]}"']&D\arrow[d,"\gamma"]\\
 C\arrow[r,"1_C"']&C.
\end{tikzcd}
\end{equation}
The occurrence with coordinates $(r,t)$ contains the placeholder
$\bigl(o_M(r),o_{o_M(r)}(t)\bigr)$.

\subsection{The finite-ordinal presentation}\label{app:ordinals}

The finite ordinals used in \cite{DP} number the elements of each
block consecutively. Put
\[
 s_0=t_0=0,\qquad
 s_r=\sum_{q=1}^{r}n_{o_M(q)}\quad(1\leq r\leq m),\qquad
 t_i=\sum_{j=1}^{i}\ell_j\quad(1\leq i\leq k).
\]
Numbering each disjoint union consecutively gives bijections
\[
 \begin{aligned}
 \theta_D:D&\longrightarrow\underline{s_m},
       &\theta_D(r,t)&=s_{r-1}+t,\\
 \theta_C:C&\longrightarrow\underline{t_k},
       &\theta_C(i,s)&=t_{i-1}+s.
 \end{aligned}
\]
Define $\delta:\underline{s_m}\to\underline{t_k}$ as the unique
function making the following square commute:
\begin{equation}\label{eq:app-ordinal}
\begin{tikzcd}[column sep=large,row sep=large]
 D\arrow[r,"\gamma"]\arrow[d,"\theta_D","\cong"']
    &C\arrow[d,"\theta_C","\cong"']\\
 \underline{s_m}\arrow[r,"\delta"']&\underline{t_k}.
\end{tikzcd}
\end{equation}
Equivalently, $\delta=\theta_C\circ\gamma\circ\theta_D^{-1}$.
For each $(r,t)\in D$, this gives
\begin{equation}\label{eq:app-coordinate}
 \delta(s_{r-1}+t)=t_{o_M(r)-1}+o_{o_M(r)}(t)
 \qquad\bigl((r,t)\in D\bigr).
\end{equation}
Every source position belongs to exactly one nonempty block, so this
formula determines $\delta$ even when some blocks are empty.

In Do\v{s}en and Petri\'c's notation, $\mu$ and $\nu_i$ correspond to
$o_M$ and $o_i$, and $\mu(\nu_1,\ldots,\nu_k)$ corresponds to $\delta$.
For $x\in\underline{s_m}$, let $r$ be the unique index with
$s_{r-1}<x\leq s_r$. Substituting $t=x-s_{r-1}$ in
\eqref{eq:app-coordinate} gives their formula
\cite[\S2.5, p.~46]{DP}. The square therefore describes the
function on finite ordinals in their generalized composition.
The \emph{skeleton} of a pattern is obtained by replacing every
placeholder by one fixed symbol and leaving the constants and
connectives unchanged. A skeleton and an occurrence map recover
the pattern by assigning the indicated placeholder to each position.
The generalized composite is therefore determined by the skeleton
of $M[\Phi_0]$ and the map $\delta$.

For example, take
\[
 M(p,q)=q\wedge(p\vee q),\qquad
 N_1(a,b)=a\wedge b,\qquad
 N_2(c,d)=d\wedge(c\vee d).
\]
The occurrence maps, listed by successive values, are
$o_M=(2,1,2)$, $o_1=(1,2)$, and $o_2=(2,1,2)$.
The source blocks have sizes $3,2,3$, while the target blocks list
$(a,b)$ and $(c,d)$. Thus
\[
\begin{array}{c|ccc|cc|ccc}
 \text{source position}&1&2&3&4&5&6&7&8\\\hline
 \text{placeholder}&d&c&d&a&b&d&c&d\\
 \delta&4&3&4&1&2&4&3&4
\end{array}
\]
For instance, source position $2$ has coordinates $(1,2)$.
The function $\gamma$ gives $\gamma(1,2)=(2,1)$, and
$\theta_C(2,1)=3$. Thus $\delta(2)=3$.

\subsection{Shared variables and the maximality argument}\label{app:shared}

The local placeholders may be assigned names in a common variable
set $Q$. Let $\rho:C\to Q$ be that assignment; different elements
of $C$ may have the same image. Replace each placeholder $c\in C$
by the variable $\rho(c)$, obtaining the renaming map
$(-)[\rho]:\FormulaSet(C)\to\FormulaSet(Q)$. Define
$\sigma:U\to\FormulaSet(Q)$ by
$\sigma(p_i)=\widetilde N_i[\rho]$.
Then the following square commutes:
\begin{equation}\label{eq:app-renaming}
\begin{tikzcd}[column sep=large,row sep=large]
 \FormulaSet(U)\arrow[r,"{(-)[\Phi_0]}"]
    \arrow[d,"{(-)[\sigma]}"']
   &\FormulaSet(C)\arrow[d,"{(-)[\rho]}"]\\
 \FormulaSet(Q)\arrow[r,"1"']&\FormulaSet(Q).
\end{tikzcd}
\end{equation}
At $p_i$, both routes give $\sigma(p_i)$. The constant and connective
clauses agree, so induction on the formula proves commutativity.
Renaming changes the labels on the occurrences. Thus the same variable
of $Q$ may occur in several inserted patterns while the pairs $(r,t)$
distinguish its occurrences.

If the literal isolated in Proposition~\ref{prop:all-contexts}
is $p$, equation~\eqref{eq:context} has source $p\wedge\bot$ and
target $p\vee\top$. Each contains one occurrence of $p$. Substituting $A$ for $p$
therefore gives one block $\Occ(A)$ at each endpoint, and the
substituted morphisms have type $A\wedge\bot\to A\vee\top$. If the isolated literal is $\neg p$,
the substitution $p\mapsto\neg A$ instead gives
$\neg\neg A\wedge\bot\to\neg\neg A\vee\top$; the
double-negation square in Proposition~\ref{prop:all-contexts}
then gives the equation for $A$.

\section{Diagrams for three defining equations}
\label{app:structural-diagrams}

\paragraph{The first-injection equation.}
For formulas $A,B,C$, equation~\eqref{eq:first-injection} is the
commutative triangle
\[
\begin{tikzcd}[column sep=large,row sep=large]
 A\wedge B
   \arrow[r,"1_A\wedge\iota_1^{B,C}"]
   \arrow[dr,"{\iota_1^{A\wedge B,C}}"']
   &A\wedge(B\vee C)\arrow[d,"d_{A,B,C}"]\\
   &(A\wedge B)\vee C.
\end{tikzcd}
\]
Both paths have source $A\wedge B$ and target $(A\wedge B)\vee C$.
The upper morphism is induced by $1_A:A\to A$ and
$\iota_1^{B,C}:B\to B\vee C$.

\paragraph{The zero identity.}
For a formula $A$, equation~\eqref{eq:zero-identity} is expressed by
the commutative diagram
\[
\begin{tikzcd}[column sep=large,row sep=large]
 A\arrow[r,"{\langle1_A,!_A\rangle}"]\arrow[d,"0_A"']
   &A\wedge\top\arrow[r,"1_A\wedge\eta_A"]
   &A\wedge(\neg A\vee A)\arrow[d,"d_{A,\neg A,A}"]\\
 A
   &\bot\vee A\arrow[l,"{[\kappa_A,1_A]}"']
   &(A\wedge\neg A)\vee A\arrow[l,"\epsilon_A\vee1_A"'].
\end{tikzcd}
\]
The path across the top, down the right side, and back along the
bottom is the composite of five morphisms in~\eqref{eq:zero-identity}.
Its first morphism pairs $1_A:A\to A$ with $!_A:A\to\top$;
its last morphism copairs $\kappa_A:\bot\to A$ with $1_A:A\to A$.
The morphisms between them use
$\eta_A:\top\to\neg A\vee A$ and
$\epsilon_A:A\wedge\neg A\to\bot$, together with the dissociativity morphism.

\paragraph{The zero law.}
For any parallel morphisms $f,g:A\to C$,
equation~\eqref{eq:zero-law} is the commutative square
\[
\begin{tikzcd}[column sep=large,row sep=large]
 A\arrow[r,"0_A"]\arrow[d,"g"']&A\arrow[d,"f"]\\
 C\arrow[r,"0_C"']&C.
\end{tikzcd}
\]
The upper-right path is $f\circ0_A:A\to C$, and the lower-left
path is $0_C\circ g:A\to C$. The morphisms $f$ and $g$ need not be equal.

\section{Corrections to the cut-elimination argument used in Boolean Coherence}
\label{app:cut-errata}

Do\v{s}en and Petri\'c's proof of Boolean Coherence
\cite[\S14.2]{DP} uses the cut-elimination theorem for distributive
lattice categories proved in \cite[\S11.2]{DP}.
Our maximality proof uses Boolean Coherence as
Theorem~\ref{thm:coherence}.
We correct two typographical errors in the rank comparisons of
that cut-elimination argument. The displayed cut reductions are unchanged.

Write $L,R$ for the left and right ranks of a cut and $r=L+R$
for its total rank. Each one-sided rank counts the nodes of the
corresponding cut-formula cluster, so $L,R\geq1$.
The induction orders $(d,r)$ lexicographically, where $d$ is the
degree of the cut formula \cite[pp.~249,253]{DP}.
In the following cases, the cut formula is unchanged.

\paragraph{Page 260, case $(w\,\gamma\,3)$.}
The inequality $r<r'$ following the displayed reduction, where $r'$
is the rank of a recursive cut, should read $r'<r$.
The two recursive cuts retain the left premise $h$.
Their right cut-formula clusters lie properly inside the original
right cluster. If their ranks are $R_f,R_g$, then
$R_f<R$ and $R_g<R$. Thus both recursive cuts have smaller total
rank: $L+R_f<r$ and $L+R_g<r$.

\paragraph{Page 262, case $(w\,\vee\,1)$.}
The phrase ``whose rank is~1'' should read ``whose right rank
is~1'', as in the analogous case on p.~261.
After eliminating the first recursive cut, the selected disjunction
in the new right premise is introduced by $\vee L$, so its
cluster has one node. The left premise is still $h$.
The new cut therefore has total rank $L+1$, not~1.
The original right cluster passes through the contraction and has
$R>1$. Hence $2\leq L+1<L+R=r$, which justifies the second
application of the induction hypothesis.

\section{The registered Lean theorems}\label{app:lean-statements}

The five declarations in the namespace \texttt{DosenPetric.Palomar}
are stated for arrow terms. We give their statements and relate them
to morphisms of $\B(P)$. The types are specified in \texttt{Challenge.lean};
\texttt{Solution.lean} imports their proofs from
\texttt{DosenPetric.PalomarMaximality}.

\subsection{Arrow terms and equations}\label{app:terms-equations}

Let $P$ be a set of propositional variables, and let
$\mathcal T_P(A,B)$ be the set of arrow terms with source
$A\in\FormulaSet(P)$ and target $B\in\FormulaSet(P)$.
Write $\mathcal E_{A,B}$ for the relation denoted by
\texttt{Derives}. It is generated by the defining equations of
the presentation, equivalence, composition, conjunction and
disjunction. Thus
\[
 \B(P)(A,B)=\mathcal T_P(A,B)/\mathcal E_{A,B}.
\]
Write $[f]$ for the class of $f\in\mathcal T_P(A,B)$.
Then $[f]=[g]$ if and only if $f\,\mathcal E_{A,B}\,g$.
The terms $\operatorname{id}_A,z_A\in\mathcal T_P(A,A)$ represent
$1_A$ and the primitive zero endomorphism, respectively; $[z_A]$ is
the morphism $0_A$ of
\eqref{eq:zero-identity}. Composition of terms is written $g\circ f$,
with $f$ applied first; in Lean this term is \texttt{Term.comp f g}.

For a set $Q$ and $\sigma:P\to\FormulaSet(Q)$, substitution gives
\[
 \mathcal T_P(A,B)\longrightarrow
 \mathcal T_Q(A[\sigma],B[\sigma]),\qquad f\longmapsto f[\sigma].
\]
Induction on an equational derivation shows that
$f\,\mathcal E_{A,B}\,g$ implies
$f[\sigma]\,\mathcal E_{A[\sigma],B[\sigma]}\,g[\sigma]$,
where the latter relation is formed over $Q$.
Consequently the following square commutes:
\[
\begin{tikzcd}[column sep=large,row sep=large]
 \mathcal T_P(A,B)\arrow[r,"{f\mapsto[f]}"]
   \arrow[d,"{f\mapsto f[\sigma]}"']
   &\B(P)(A,B)\arrow[d,"{[f]\mapsto F_\sigma([f])}"]\\
 \mathcal T_Q(A[\sigma],B[\sigma])
   \arrow[r,"{u\mapsto[u]}"']
   &\B(Q)(A[\sigma],B[\sigma]).
\end{tikzcd}
\]

A \texttt{PrimitiveCongruence P} is a family of equivalence relations
$R_{A,B}$ on $\mathcal T_P(A,B)$ that contains every defining equation
and is preserved by composition, conjunction and disjunction.
Explicitly, composition requires
\[
 \begin{gathered}
 f\,R_{A,B}\,f',\quad g\,R_{B,C}\,g'
 \quad\Longrightarrow\quad
 g\circ f\,R_{A,C}\,g'\circ f',
 \end{gathered}
\]
and, for $\star\in\{\wedge,\vee\}$, the binary operations require
\[
 f\,R_{A,B}\,f',\quad g\,R_{C,D}\,g'
 \quad\Longrightarrow\quad
 f\star g\,R_{A\star C,B\star D}\,f'\star g'.
\]
The property \texttt{R.Closed} is preservation under every
$\sigma:P\to\FormulaSet(P)$:
\[
 f\,R_{A,B}\,g
 \quad\Longrightarrow\quad
 f[\sigma]\,R_{A[\sigma],B[\sigma]}\,g[\sigma].
\]
The property \texttt{R.Total} means
\[
 \text{for every }A,B\text{ and every }f,g\in\mathcal T_P(A,B),
 \qquad f\,R_{A,B}\,g.
\]
These quantifiers do not assert that $\mathcal T_P(A,B)$ is nonempty.

\subsection{Reading a Lean declaration}\label{app:reading-lean}

The Lean identifiers below refer to the registered source archive
\cite{palomar-2026-10-30-000014-v3}.
Lean represents a mathematical set by a type and a proposition by
an expression of type \texttt{Prop}. A proof of a proposition is a
term having that proposition as its type. Thus
\texttt{f : Term P A B} declares an arrow term with source $A$
and target $B$, whereas \texttt{h : Derives f g} declares a proof
that the two terms are related by the presentation's equations.
The notation \texttt{Hom A B} denotes the quotient morphisms, rather
than the arrow terms.\footnote{\raggedright For the general syntax, see
\emph{Theorem Proving in Lean 4}, Chapters 2--3 and 5,
\url{https://lean-lang.org/theorem_proving_in_lean4/}.}

For example, the statement of \href{https://github.com/flengyel/Dosen-Petric-Maximality/blob/dosen-petric-v25-palomar-v3/DosenPetric/PalomarMaximality.lean\#L17}{\texttt{Palomar.maximality}},
writing \texttt{Not} for the notation $\neg$, is
\begin{quote}
\begin{verbatim}
theorem maximality [Nonempty P]
    (R : PrimitiveCongruence P) (closed : R.Closed)
    {A B : Formula P} (f g : Term P A B)
    (hne : Not (Derives f g)) (hfg : R.Rel f g) :
    R.Total
\end{verbatim}
\end{quote}
Here $P$ is a parameter declared in the surrounding scope as
\texttt{P : Type u}. This represents the variable set at universe
level \texttt{u}; the theorem is uniform in that level.
The bracketed \texttt{[Nonempty P]} is a hypothesis that $P$ is
nonempty. The braces in \texttt{\{A B : Formula P\}} mark
universally quantified arguments that Lean infers from the types
of \texttt{f} and \texttt{g}.
The names \texttt{closed}, \texttt{hne} and \texttt{hfg} name
proofs of the three displayed hypotheses; the expression after
the last colon is the conclusion.

Dot notation has two uses here. In
\href{https://github.com/flengyel/Dosen-Petric-Maximality/blob/dosen-petric-v25-palomar-v3/DosenPetric/PalomarDefinitions.lean\#L14}{\texttt{Palomar.PrimitiveCongruence}}, it gives the namespace of a
definition. In \texttt{R.Rel}, it selects the relation belonging
to $R$. The expression \texttt{R.Closed} applies the defined
property \href{https://github.com/flengyel/Dosen-Petric-Maximality/blob/dosen-petric-v25-palomar-v3/DosenPetric/PalomarDefinitions.lean\#L31}{\texttt{PrimitiveCongruence.Closed}} to $R$.
The assertion \texttt{R.Total} expands to
\begin{quote}
\begin{verbatim}
forall {C D} (k l : Term P C D), R.Rel k l
\end{verbatim}
\end{quote}
which is the universal statement in the preceding subsection.
Consequently the theorem takes the given relation and hypotheses
to a proof that every pair of parallel terms is related.

In a proof introduced by \texttt{by}, the commands construct this
proof term. The command \texttt{intro} introduces a universally
quantified variable or a hypothesis; \texttt{obtain} extracts
witnesses and their properties from an existence statement;
\texttt{rw} replaces expressions using equations; and
\texttt{exact} supplies a proof of the statement remaining to
be proved. Lean's kernel checks that the resulting term has the theorem's displayed type.
In particular, \texttt{Term.comp f g} applies \texttt{f} first,
so it denotes $[g]\circ[f]$ on quotient morphisms.

The proof of \href{https://github.com/flengyel/Dosen-Petric-Maximality/blob/dosen-petric-v25-palomar-v3/DosenPetric/PalomarMaximality.lean\#L27}{\texttt{Palomar.short\_maximality}} can be read as
the following sequence of mathematical implications.
\begin{center}
\small
\begin{tabular}{>{\raggedright\arraybackslash}p{0.37\linewidth}>{\raggedright\arraybackslash}p{0.55\linewidth}}
\textit{Lean declaration} & \textit{Mathematical step}\\[3pt]
\href{https://github.com/flengyel/Dosen-Petric-Maximality/blob/dosen-petric-v25-palomar-v3/DosenPetric/PalomarCorrespondence.lean\#L121}{\texttt{Palomar.PrimitiveCongruence.\allowbreak{}raw\_distinct\_iff}}
 & $\neg\mathcal E(f,g)$ means $[f]\ne[g]$.\\[5pt]
\href{https://github.com/flengyel/Dosen-Petric-Maximality/blob/dosen-petric-v25-palomar-v3/DosenPetric/Isolation.lean\#L219}{\texttt{Isolation.\allowbreak{}normalized\_literal\_context}}
 & Lemma~\ref{lem:isolation} and compatibility with composition give a literal $\ell$ and $(K_\ell)$
   modulo the induced congruence.\\[5pt]
\href{https://github.com/flengyel/Dosen-Petric-Maximality/blob/dosen-petric-v25-palomar-v3/DosenPetric/Substitution.lean\#L90}{\texttt{ContextualCollapse.\allowbreak{}all\_contexts}}
 & Proposition~\ref{prop:all-contexts} gives $(K_A)$ for every formula $A$.\\[5pt]
\href{https://github.com/flengyel/Dosen-Petric-Maximality/blob/dosen-petric-v25-palomar-v3/DosenPetric/PresentationShortProjection.lean\#L86}{\texttt{Presentation.ShortProof.\allowbreak{}total}}
 & Lemma~\ref{lem:collapse} identifies arbitrary parallel morphisms by
   the $q_A,e_A$ argument.
\end{tabular}
\end{center}
In this proof \texttt{R.fragment} is the induced relation on
quotient morphisms, with its composition, pairing and copairing
laws. The argument \texttt{interpretation P} supplies the
occurrence interpretation and its proved faithfulness;
\texttt{normalization P} supplies the normal-form isomorphisms.
The arguments \texttt{substitution P} and
\texttt{doubleNegation P} supply the morphism substitution and
double-negation maps used in Proposition~\ref{prop:all-contexts}.

\subsection{The five declarations}

\paragraph{\texttt{maximality} and \texttt{short\_maximality}.}
Each declaration has the following statement. Suppose that $P$ is
nonempty, $R$ is as above and \texttt{R.Closed} holds. If, for some
$A,B$ and $f,g\in\mathcal T_P(A,B)$,
\[
 f\,R_{A,B}\,g
 \qquad\text{and}\qquad
 \neg\bigl(f\,\mathcal E_{A,B}\,g\bigr),
\]
then \texttt{R.Total} holds.

The declaration \texttt{maximality} uses the argument with the
symmetric morphism $h_A$ in Section~\ref{sec:formalization}.
The declaration \texttt{short\_maximality} uses the isomorphisms
$q_A,e_A$ in Section~\ref{sec:collapse}.
Both use the occurrence selection of Lemma~\ref{lem:isolation}
and substitution as in Proposition~\ref{prop:all-contexts}.
The relation used in these arguments need only be preserved by
composition, pairing and copairing, together with uniform substitution.
Pairing and copairing preserve the relation induced by $R$ because
they are defined using the diagonal, codiagonal and the two binary
operations on morphisms.

\paragraph{\texttt{congruence\_dichotomy}.}
For nonempty $P$ and a family $R$ satisfying \texttt{R.Closed},
either
\[
 \text{for every }A,B\text{ and }f,g\in\mathcal T_P(A,B),\qquad
 f\,R_{A,B}\,g\ \Longleftrightarrow\ f\,\mathcal E_{A,B}\,g,
\]
or \texttt{R.Total} holds. This is the equality-or-totality formulation
of Theorem~\ref{thm:main}.
Indeed, induction on \texttt{Derives} gives $\mathcal E_{A,B}\subseteq
R_{A,B}$. If the reverse inclusion fails for some pair, that pair
satisfies the hypotheses of \texttt{maximality}.

\paragraph{\texttt{identity\_ne\_zero}.}
For every $p\in P$,
\[
 \neg\bigl(\operatorname{id}_p\,\mathcal E_{p,p}\,z_p\bigr).
\]
Equivalently, $1_p\ne0_p$ in $\B(P)$.
The occurrence relation of $1_p$ is the diagonal on the singleton
$\Occ(p)$, whereas the occurrence relation of $0_p$ is empty.
This is the distinction used in Section~\ref{sec:occurrences}.
The declaration takes $p$ as an argument and makes no separate
assumption that $P$ is nonempty.

\paragraph{\texttt{no\_top\_to\_bot}.}
For every set $P$, including the empty set,
\[
 \mathcal T_P(\top,\bot)=\varnothing.
\]
Since every morphism is represented by a term, this is equivalent to
$\B(P)(\top,\bot)=\varnothing$, as stated in
Section~\ref{sec:formalization}.
Let $\tau:\FormulaSet(P)\to\{0,1\}$ be evaluation by the Boolean
operations with $\tau(p)=1$ for every variable $p$.
Induction on $f\in\mathcal T_P(A,B)$ gives $\tau(A)\leq\tau(B)$
in the order $0<1$. A term of type
$\top\to\bot$ would give $1\leq0$, which is impossible.

\subsection{Congruences on terms and on morphisms}

For a family $R$ satisfying the conditions above,
$\mathcal E_{A,B}\subseteq R_{A,B}$ implies that
\[
 [f]\,\overline R_{A,B}\,[g]
 \quad\Longleftrightarrow\quad f\,R_{A,B}\,g
\]
defines an equivalence relation on $\B(P)(A,B)$, independent of
representatives. For example, if
$f\,\mathcal E_{A,B}\,f'$ and $g\,\mathcal E_{A,B}\,g'$, then
equivalence and $\mathcal E_{A,B}\subseteq R_{A,B}$ give
$f\,R_{A,B}\,g$ if and only if $f'\,R_{A,B}\,g'$.
The relations $\overline R_{A,B}$ are preserved by composition,
pairing and copairing. If \texttt{R.Closed} holds, the substitution
square in Appendix~\ref{app:terms-equations}, with $Q=P$, shows
that they are preserved by $F_\sigma$ for
every $\sigma:P\to\FormulaSet(P)$.

The relation $\overline R$ is equality on morphisms if and only if
$R=\mathcal E$, and totality of $\overline R$ is equivalent to
\texttt{R.Total}. If $P$ is nonempty and \texttt{R.Closed} holds, the dichotomy
therefore makes $\overline R$ either equality or the total relation. Both are preserved
by every operation on morphisms, including negation. Thus $\overline R$ is a
substitution congruence in the sense of Section~\ref{sec:category}.

Conversely, let $P$ be nonempty and let $\theta$ be a substitution
congruence on $\B(P)$ in the sense of Section~\ref{sec:category}.
Define
\[
 f\,R^\theta_{A,B}\,g
 \quad\Longleftrightarrow\quad [f]\,\theta_{A,B}\,[g].
\]
The family $R^\theta$ satisfies the defining conditions of
\texttt{PrimitiveCongruence P} and \texttt{Closed}.
Distinct $\theta$-equivalent morphisms have representatives that are
$R^\theta$-related but not $\mathcal E$-related. Either maximality
declaration then makes $R^\theta$ total, and hence makes $\theta$ total.
The quotient $\B(P)/\theta$ is therefore a preorder, as stated in
Theorem~\ref{thm:main}.
The correspondence is proved in
\texttt{DosenPetric.\allowbreak{}PalomarCorrespondence}.

The theorem \texttt{Presentation.\allowbreak{}coherence} proves
Theorem~\ref{thm:coherence}, and
\texttt{PresentationSemantics.\allowbreak{}model\_corollary} proves
Corollary~\ref{cor:models}.
Both names have prefix \texttt{DosenPetric.}

In the namespace \texttt{DosenPetric.\allowbreak{}ContextualCollapse},
\texttt{second\_block} proves
$j_A\circ h_A\equiv j_A\circ\pi_2^{A,A}$, and
\texttt{fold\_top} proves
$h_A\equiv\pi_1^{A,A}$, as calculated in
Section~\ref{sec:formalization}.
The intermediate equality making $A\vee\top$ terminal is
\texttt{terminal\_retraction}.

\subsection{The proof of coherence}\label{app:coherence-proof}

Let $\mathcal R(P)$ denote the presentation on negation normal
forms described in Section~\ref{sec:occurrences}. We use the same
notation for its products, coproducts, units, associativity and
symmetry maps, mix and dissociativity. Write $G$ also for its occurrence
interpretation; inclusion into $\B(P)$ preserves each tagged
relation. A term of $\mathcal R(P)$
is \emph{positive} if it contains no $\eta_p$ or $\epsilon_p$
constructor. A morphism is positive when it has such a representative.

The union equations used in the argument below have derivations
from the defining equations of $\mathcal R(P)$.
Naturality of the diagonal, mix and codiagonal gives compatibility
with composition; symmetry and associativity give commutativity
and associativity. The equation
\[
 [1_A,1_A]\circ m_{A,A}\circ\Delta_A=1_A
\]
gives idempotence, and the zero law and the product, coproduct
and dissociativity equations give $f\cup0=f$.
These derivations do not use faithfulness.\footnote{\raggedright Lean:
\href{https://github.com/flengyel/Dosen-Petric-Maximality/blob/dosen-petric-v25-palomar-v3/DosenPetric/ReducedUnion.lean\#L80}{\texttt{Reduced.\allowbreak{}union\_pre}},
\href{https://github.com/flengyel/Dosen-Petric-Maximality/blob/dosen-petric-v25-palomar-v3/DosenPetric/ReducedUnion.lean\#L85}{\texttt{Reduced.\allowbreak{}union\_post}},
\href{https://github.com/flengyel/Dosen-Petric-Maximality/blob/dosen-petric-v25-palomar-v3/DosenPetric/ReducedUnion.lean\#L90}{\texttt{Reduced.\allowbreak{}union\_comm}},
\href{https://github.com/flengyel/Dosen-Petric-Maximality/blob/dosen-petric-v25-palomar-v3/DosenPetric/ReducedUnion.lean\#L131}{\texttt{Reduced.\allowbreak{}union\_assoc}},
\href{https://github.com/flengyel/Dosen-Petric-Maximality/blob/dosen-petric-v25-palomar-v3/DosenPetric/ReducedUnion.lean\#L97}{\texttt{Reduced.\allowbreak{}union\_idem}}, and
\href{https://github.com/flengyel/Dosen-Petric-Maximality/blob/dosen-petric-v25-palomar-v3/DosenPetric/ReducedUnion.lean\#L145}{\texttt{Reduced.\allowbreak{}union\_id\_zero}}.}

The two canonical comparison morphisms
\[
 \begin{aligned}
 c^\wedge_{X,U,V}
 &=[1_X\wedge\iota_1^{U,V},1_X\wedge\iota_2^{U,V}]
 :(X\wedge U)\vee(X\wedge V)\longrightarrow X\wedge(U\vee V),\\
 c^\vee_{U,V,X}
 &=\langle\pi_1^{U,V}\vee1_X,\pi_2^{U,V}\vee1_X\rangle
 :(U\wedge V)\vee X\longrightarrow(U\vee X)\wedge(V\vee X)
 \end{aligned}
\]
have, respectively, a section $s_{X,U,V}$ and a retraction
$t_{U,V,X}$:
\[
 c^\wedge_{X,U,V}\circ s_{X,U,V}=1_{X\wedge(U\vee V)},\qquad
 t_{U,V,X}\circ c^\vee_{U,V,X}=1_{(U\wedge V)\vee X}.
\]
Appendix~\ref{app:explicit-inverses} defines these maps by
composites, products and coproducts of identities, associativity,
symmetry and dissociativity maps, diagonals and codiagonals.
Their two equations follow from the
product, coproduct, dissociativity and mix equations.\footnote{\raggedright Lean:
\href{https://github.com/flengyel/Dosen-Petric-Maximality/blob/dosen-petric-v25-palomar-v3/DosenPetric/GentzenDistributiveEpi.lean\#L57}{\texttt{Reduced.Gentzen.DistributiveEpi.\allowbreak{}section\_comparison}}
and
\href{https://github.com/flengyel/Dosen-Petric-Maximality/blob/dosen-petric-v25-palomar-v3/DosenPetric/GentzenDistributiveMono.lean\#L54}{\texttt{Reduced.Gentzen.DistributiveMono.\allowbreak{}comparison\_retraction}}.}

Iterating the first split comparison gives, for every formula $A$
of $\mathcal R(P)$, a disjunctive normal-form expansion $D_A$
and positive morphisms
\[
 i_A:A\longrightarrow D_A,\qquad r_A:D_A\longrightarrow A,
 \qquad r_A\circ i_A=1_A.
\]
Here $D_A$ is a coproduct of conjunctions of literals and constants.
The second split comparison gives a conjunctive expansion $C_A$,
a product of disjunctions of literals and constants, with positive
morphisms $A\to C_A\to A$ whose composite is $1_A$.
Induction on $A$ constructs these expansions: conjunctions of
disjunctive expansions use $s$, disjunctions of conjunctive
expansions use $c^\vee$, and the other connective combines
expansions by its bifunctor.
Only the displayed retraction equations are needed.\footnote{\raggedright Lean:
\href{https://github.com/flengyel/Dosen-Petric-Maximality/blob/dosen-petric-v25-palomar-v3/DosenPetric/GentzenNormalExpansion.lean\#L45}{\texttt{Reduced.Gentzen.\allowbreak{}dnf\_expansion}} and
\href{https://github.com/flengyel/Dosen-Petric-Maximality/blob/dosen-petric-v25-palomar-v3/DosenPetric/GentzenNormalExpansion.lean\#L62}{\texttt{Reduced.Gentzen.\allowbreak{}cnf\_expansion}}.}

Consider a positive morphism from a conjunction $U$ of literals
and constants to a disjunction $V$ of literals and constants.
At occurrences $x\in\Occ(U)$ and $y\in\Occ(V)$ of the same
literal $\ell$, the successive product projections and coproduct
injections give
\[
 p_x:U\longrightarrow\ell,\qquad
 j_y:\ell\longrightarrow V.
\]
The relation of $j_y\circ p_x$ is the singleton $\{(x,y)\}$.
Represent the given morphism by a Gentzen derivation and eliminate
its cuts. Replace weakening on compound formulas by weakening on
literals and constants. The normal-form theorem then expresses
the original morphism as a finite union of these paths; an empty
union uses the zero in its inhabited hom-set. The relation of that union records
exactly its paths. Idempotence and commutativity make repetitions
and their order irrelevant, so equal occurrence relations give
equal positive morphisms at these endpoints.\footnote{\raggedright Lean:
\href{https://github.com/flengyel/Dosen-Petric-Maximality/blob/dosen-petric-v25-palomar-v3/DosenPetric/GentzenMixNormalForm.lean\#L37}{\texttt{Reduced.Gentzen.G.\allowbreak{}mix\_normal\_form}},
\href{https://github.com/flengyel/Dosen-Petric-Maximality/blob/dosen-petric-v25-palomar-v3/DosenPetric/GentzenPureNormalForms.lean\#L84}{\texttt{Reduced.Gentzen.HasPureNormalForm.\allowbreak{}faithful}}, and
\href{https://github.com/flengyel/Dosen-Petric-Maximality/blob/dosen-petric-v25-palomar-v3/DosenPetric/GentzenPurePositiveCoherence.lean\#L14}{\texttt{Reduced.Gentzen.\allowbreak{}pure\_positive\_coherence}}.}
The coproduct and product equations extend this assertion to
morphisms from a disjunctive expansion to a conjunctive expansion.

For positive $f,g:A\to B$, choose a disjunctive expansion of $A$
and a conjunctive expansion of $B$. Write their maps as
$i_A:A\to D_A$, $r_A:D_A\to A$ and
$i_B:B\to C_B$, $r_B:C_B\to B$, with
$r_A\circ i_A=1_A$ and $r_B\circ i_B=1_B$.
If $Gf=Gg$, the morphisms
$i_B\circ f\circ r_A,i_B\circ g\circ r_A:D_A\to C_B$
have equal relations and hence are equal. Precomposing with
$i_A$ and postcomposing with $r_B$ gives $f=g$.
This proves positive coherence at arbitrary endpoints.\footnote{\raggedright Lean:
\href{https://github.com/flengyel/Dosen-Petric-Maximality/blob/dosen-petric-v25-palomar-v3/DosenPetric/GentzenPositiveCoherence.lean\#L14}{\texttt{Reduced.Gentzen.\allowbreak{}represented\_coherence}} and
\href{https://github.com/flengyel/Dosen-Petric-Maximality/blob/dosen-petric-v25-palomar-v3/DosenPetric/GentzenPositiveCoherence.lean\#L31}{\texttt{Reduced.Gentzen.\allowbreak{}positive\_coherence}}.}

For arbitrary terms $f,g:A\to B$ of $\mathcal R(P)$, there are
formulas $A',B'$, morphisms
$u:A\to A'$, $v:B'\to B$ and positive terms
$f_0,g_0:A'\to B'$ such that, on quotient morphisms,
\[
 [f]=v\circ[f_0]\circ u,\qquad [g]=v\circ[g_0]\circ u.
\]
There are also injections
$a:\Occ(A)\to\Occ(A')$ and $b:\Occ(B)\to\Occ(B')$ with
\[
 G[f_0]=(a\times b)(G[f]),\qquad
 G[g_0]=(a\times b)(G[g]).
\]
The right sides are direct images of relations.

To construct these data, factor $f\wedge g:A\wedge A\to B\wedge B$
by induction on the number of $\eta_p$ and $\epsilon_p$
constructors. Extracting an $\eta_p$ adds $\neg p\vee p$ as a
conjunct at the source; extracting an $\epsilon_p$ adds
$p\wedge\neg p$ as a disjunct at the target. Each step removes
one such constructor from the middle term, leaving a positive
middle term when the induction terminates. Extend the diagonal
$\Delta_A:A\to A\wedge A$ and the projections
$\pi_i^{B,B}:B\wedge B\to B$, for $i=1,2$, by identities on
the added formulas. Composing with these extensions gives two
positive terms with source $A'$ and target $B'$; naturality of
the extracted prefixes and suffixes gives the common morphisms
$u$ and $v$. Finally, precompose and postcompose these terms
with endomorphisms of $A'$ and $B'$ formed from $1_A$ and
$1_B$, respectively, and zero endomorphisms on the added
conjunction factors and disjunction summands.
These endomorphisms preserve the factorizations and give the
displayed direct-image relations.

Thus equality of $G[f]$ and $G[g]$ gives equality of the middle relations.
Positive coherence equates $[f_0]$ and $[g_0]$, and the
factorizations give $[f]=[g]$.\footnote{\raggedright Lean:
\href{https://github.com/flengyel/Dosen-Petric-Maximality/blob/dosen-petric-v25-palomar-v3/DosenPetric/CommonAtomicFactorization.lean\#L51}{\texttt{Reduced.\allowbreak{}commonAtomicFactor}} and
\href{https://github.com/flengyel/Dosen-Petric-Maximality/blob/dosen-petric-v25-palomar-v3/DosenPetric/CommonAtomicFactorization.lean\#L123}{\texttt{Reduced.\allowbreak{}faithful\_of\_positive\_coherence}}.}
The equivalence described in Section~\ref{sec:occurrences}
preserves the occurrence interpretation and therefore transfers
faithfulness from $\mathcal R(P)$ to $\B(P)$.\footnote{\raggedright Lean:
\href{https://github.com/flengyel/Dosen-Petric-Maximality/blob/dosen-petric-v25-palomar-v3/DosenPetric/ReducedCoherence.lean\#L17}{\texttt{Reduced.\allowbreak{}boolean\_coherence}} and
\href{https://github.com/flengyel/Dosen-Petric-Maximality/blob/dosen-petric-v25-palomar-v3/DosenPetric/PresentationCoherence.lean\#L13}{\texttt{Presentation.\allowbreak{}coherence}}.}

\section{Inverses of \texorpdfstring{$q_A$ and $e_A$}{qA and eA}}
\label{app:explicit-inverses}

We construct inverses for the morphisms in
Lemma~\ref{lem:structural-isos}. The inverse identities below
use Boolean Coherence, whereas the split-comparison identities
required to prove it follow directly from the presentation equations.
Throughout this appendix,
\[
 \alpha^\star_{X,Y,Z}:(X\star Y)\star Z\longrightarrow
          X\star(Y\star Z),\qquad
 \beta^\star_{X,Y}:X\star Y\longrightarrow Y\star X
 \quad(\star\in\{\wedge,\vee\})
\]
are associativity and symmetry.

\subsection{The inverse of \texorpdfstring{$q_A$}{qA}}

Put $S_A=(A\wedge\bot)\vee\top$. Define $u_A:A\to S_A$ by
the path
\[
 A\xrightarrow{\langle1_A,!_A\rangle}A\wedge\top
  \xrightarrow{1_A\wedge\iota_2^{\bot,\top}}A\wedge(\bot\vee\top)
  \xrightarrow{d_{A,\bot,\top}}S_A,
\]
and put
\[
 \bar q_A=[u_A,\iota_2^{A\wedge\bot,\top}]:
      A\vee\top\longrightarrow S_A.
\]
The occurrence relation of $\bar q_A$ is the graph sending
$(1,x)\in\Occ(A\vee\top)$ to
$(1,(1,x))\in\Occ(S_A)$; $Gq_A$ is the inverse graph.
Boolean Coherence therefore gives
\[
 \bar q_A\circ q_A=1_{S_A},\qquad
 q_A\circ\bar q_A=1_{A\vee\top}.
\]
Thus $\bar q_A=q_A^{-1}$.\footnote{\raggedright Lean:
\href{https://github.com/flengyel/Dosen-Petric-Maximality/blob/dosen-petric-v25-palomar-v3/DosenPetric/PresentationShortQ.lean\#L17}{\texttt{Presentation.ShortProof.\allowbreak{}qInverseTerm}},
\href{https://github.com/flengyel/Dosen-Petric-Maximality/blob/dosen-petric-v25-palomar-v3/DosenPetric/PresentationShortQ.lean\#L32}{\texttt{Presentation.ShortProof.\allowbreak{}q\_qInverse}}, and
\href{https://github.com/flengyel/Dosen-Petric-Maximality/blob/dosen-petric-v25-palomar-v3/DosenPetric/PresentationShortQ.lean\#L46}{\texttt{Presentation.ShortProof.\allowbreak{}qInverse\_q}}.}

\subsection{The inverse of \texorpdfstring{$e_A$}{eA}}

Define $c_{X,Y;U,V}:(X\wedge Y)\wedge(U\vee V)
\to(X\wedge U)\vee(Y\wedge V)$ by the path
\[
\begin{tikzcd}[row sep=large]
 (X\wedge Y)\wedge(U\vee V)
   \arrow[d,"{1_{X\wedge Y}\wedge\beta^\vee_{U,V}}"]\\
 (X\wedge Y)\wedge(V\vee U)
   \arrow[d,"{\alpha^\wedge_{X,Y,V\vee U}}"]\\
 X\wedge(Y\wedge(V\vee U))
   \arrow[d,"{1_X\wedge d_{Y,V,U}}"]\\
 X\wedge((Y\wedge V)\vee U)
   \arrow[d,"{1_X\wedge\beta^\vee_{Y\wedge V,U}}"]\\
 X\wedge(U\vee(Y\wedge V))
   \arrow[d,"{d_{X,U,Y\wedge V}}"]\\
 (X\wedge U)\vee(Y\wedge V).
\end{tikzcd}
\]
Set
\[
 s_{X,U,V}=c_{X,X;U,V}\circ(\Delta_X\wedge1_{U\vee V}):
 X\wedge(U\vee V)\longrightarrow(X\wedge U)\vee(X\wedge V).
\]
With $T_A=A\wedge\top$, define
\[
 \bar e_A=[1_{A\wedge A},0_{T_A,A\wedge A}]
                 \circ s_{A,A,\top}:
 A\wedge(A\vee\top)\longrightarrow A\wedge A.
\]
The composite $\Delta_A\circ\pi_1^{A,\top}:T_A\to A\wedge A$
supplies the witness defining this zero.
The relation of $s_{A,A,\top}$ copies each occurrence of its
first $A$ into both target summands and sends the occurrences
of the second $A$ to the second factor of the first summand.
Postcomposition with $[1_{A\wedge A},0_{T_A,A\wedge A}]$
removes the second summand. Hence $G\bar e_A$ sends
$(1,x)$ to $(1,x)$ and $(2,(1,x))$ to $(2,x)$.
These assignments are inverse to $Ge_A$. Coherence gives
\[
 \bar e_A\circ e_A=1_{A\wedge A},\qquad
 e_A\circ\bar e_A=1_{A\wedge(A\vee\top)}.
\]
Thus $\bar e_A=e_A^{-1}$.\footnote{\raggedright Lean:
\href{https://github.com/flengyel/Dosen-Petric-Maximality/blob/dosen-petric-v25-palomar-v3/DosenPetric/PresentationShortE.lean\#L13}{\texttt{Presentation.ShortProof.\allowbreak{}sourceConjTerm}},
\href{https://github.com/flengyel/Dosen-Petric-Maximality/blob/dosen-petric-v25-palomar-v3/DosenPetric/PresentationShortE.lean\#L21}{\texttt{Presentation.ShortProof.\allowbreak{}sectionTerm}},
\href{https://github.com/flengyel/Dosen-Petric-Maximality/blob/dosen-petric-v25-palomar-v3/DosenPetric/PresentationShortE.lean\#L30}{\texttt{Presentation.ShortProof.\allowbreak{}eInverseTerm}},
\href{https://github.com/flengyel/Dosen-Petric-Maximality/blob/dosen-petric-v25-palomar-v3/DosenPetric/PresentationShortE.lean\#L44}{\texttt{Presentation.ShortProof.\allowbreak{}e\_eInverse}}, and
\href{https://github.com/flengyel/Dosen-Petric-Maximality/blob/dosen-petric-v25-palomar-v3/DosenPetric/PresentationShortE.lean\#L50}{\texttt{Presentation.ShortProof.\allowbreak{}eInverse\_e}}.}
Both inverse constructions apply to every formula $A$, including
formulas with empty occurrence sets.

\subsection{The split comparisons used in coherence}

The map $s_{X,U,V}$ is also a section of $c^\wedge_{X,U,V}$ in
Appendix~\ref{app:coherence-proof}. To construct the retraction
of $c^\vee_{U,V,X}$, first put
\[
 v_{U,X;V}=(\beta^\wedge_{V,U}\vee1_X)\circ d_{V,U,X}
                \circ\beta^\wedge_{U\vee X,V}:
 (U\vee X)\wedge V\longrightarrow(U\wedge V)\vee X.
\]
For another formula $Y$, put
\[
 \begin{aligned}
 b_{U,X;V,Y}
 &=\alpha^\vee_{U\wedge V,X,Y}\circ(v_{U,X;V}\vee1_Y)
                      \circ d_{U\vee X,V,Y}\\
 &:(U\vee X)\wedge(V\vee Y)
                   \longrightarrow(U\wedge V)\vee(X\vee Y).
 \end{aligned}
\]
Then
\[
 t_{U,V,X}=(1_{U\wedge V}\vee[1_X,1_X])\circ b_{U,X;V,X}:
 (U\vee X)\wedge(V\vee X)\longrightarrow(U\wedge V)\vee X.
\]
Expanding the definitions and using the presentation equations
gives
\[
 \begin{aligned}
 \pi_i^{X,U\vee V}\circ c^\wedge_{X,U,V}\circ s_{X,U,V}
     &=\pi_i^{X,U\vee V}\quad(i=1,2),\\
 t_{U,V,X}\circ c^\vee_{U,V,X}\circ\iota_i^{U\wedge V,X}
     &=\iota_i^{U\wedge V,X}\quad(i=1,2).
 \end{aligned}
\]
The product and coproduct equations give the two split identities.
For the first projection in the first line and the second injection
in the second line, the calculation uses
$[1_Z,1_Z]\circ m_{Z,Z}\circ\Delta_Z=1_Z$.
These calculations use only the product, coproduct, associativity,
symmetry, dissociativity and mix equations, so the same expressions
and identities hold in $\mathcal R(P)$.

\section*{Use of AI tools}
OpenAI's ChatGPT assisted with the development of the proof, preparation
of the manuscript, and checks of the arguments and references. OpenAI Codex
was used to develop the Lean formalization and to run its verification
checks. The author takes responsibility for the mathematical claims and
the final text.

\clearpage
\begingroup
\small
\setlength{\bibsep}{2pt}
\bibliographystyle{unsrtnat}
\bibliography{references_v27,palomar_v27}

@misc{palomar-2026-10-30-000014-v3,
  author = {{Florian Lengyel}},
  title = {{Maximality of the Došen–Petrić Boolean Category}},
  year = {2026},
  howpublished = {Palomar, PALOMAR-2026-10-30-000014 v3},
  url = {https://palomar-registry.org/entry?id=PALOMAR-2026-10-30-000014&version=3},
}

@book{DP,
  author = {Kosta Do{\v{s}}en and Zoran Petri{\'c}},
  title = {{Proof-Theoretical Coherence}},
  series = {Studies in Logic},
  volume = {1},
  publisher = {King's College Publications},
  address = {London},
  year = {2004},
  note = {Revised author version, September 2007},
  url = {https://www.mi.sanu.ac.rs/novi_sajt/biography/kosta/coh.pdf},
}

@article{DPCartesian,
  author = {Kosta Do{\v{s}}en and Zoran Petri{\'c}},
  title = {{The maximality of cartesian categories}},
  journal = {Mathematical Logic Quarterly},
  volume = {47},
  year = {2001},
  pages = {137--144},
  url = {https://arxiv.org/abs/math/9911059},
}

@misc{DPBicartesian,
  author = {Kosta Do{\v{s}}en and Zoran Petri{\'c}},
  title = {{Bicartesian coherence revisited}},
  howpublished = {In Zoran Ognjanovi{\'c} (ed.), \emph{Logic in Computer Science}, Zbornik radova \textbf{12(20)}, pages 5--34},
  year = {2009},
  note = {Author version: arXiv:0711.4961v3, December 8, 2008},
  url = {https://arxiv.org/abs/0711.4961},
}

@article{DPIsomorphic,
  author = {Kosta Do{\v{s}}en and Zoran Petri{\'c}},
  title = {{Isomorphic formulae in classical propositional logic}},
  journal = {Mathematical Logic Quarterly},
  volume = {58},
  year = {2012},
  pages = {5--17},
  note = {Author version: \url{https://arxiv.org/abs/0912.1271}, arXiv:0912.1271v3, October 4, 2010},
  url = {https://doi.org/10.1002/malq.201020020},
}

@article{DPSplit,
  author = {Kosta Do{\v{s}}en and Zoran Petri{\'c}},
  title = {{Syntax for split preorders}},
  journal = {Annals of Pure and Applied Logic},
  volume = {164},
  year = {2013},
  pages = {443--481},
  note = {Author version: arXiv:0902.0742v7, October 18, 2012},
  url = {https://arxiv.org/abs/0902.0742},
}

@inproceedings{LSFree,
  author = {Fran{\c{c}}ois Lamarche and Lutz Stra{\ss}burger},
  title = {{Constructing free Boolean categories}},
  booktitle = {Proceedings of the 20th Annual IEEE Symposium on Logic in Computer Science (LICS 2005)},
  publisher = {IEEE Computer Society},
  year = {2005},
  pages = {209--218},
  url = {https://doi.org/10.1109/LICS.2005.13},
  note = {\href{https://www.lix.polytechnique.fr/~lutz/papers/FreeBool.pdf}{Author's proceedings version}},
}

@article{FP,
  author = {Carsten F{\"u}hrmann and David Pym},
  title = {{Order-enriched categorical models of the classical sequent calculus}},
  journal = {Journal of Pure and Applied Algebra},
  volume = {204},
  year = {2006},
  pages = {21--78},
  url = {https://doi.org/10.1016/j.jpaa.2005.03.016},
  note = {\href{https://carstenfuehrmann.org/papers/order-enriched-models-for-classical-sequent-calculus.pdf}{Author version}},
}

@article{Strassburger,
  author = {Lutz Stra{\ss}burger},
  title = {{On the axiomatisation of Boolean categories with and without medial}},
  journal = {Theory and Applications of Categories},
  volume = {18},
  year = {2007},
  pages = {536--601},
  url = {https://www.tac.mta.ca/tac/volumes/18/18/18-18abs.html},
}

@misc{Formalization,
  author = {Florian Lengyel},
  title = {{Do{\v{s}}en--Petri{\'c} maximality in Lean}},
  howpublished = {Lean 4 source development},
  year = {2026},
  note = {Release \nolinkurl{dosen-petric-v25-palomar-v3}},
  url = {https://github.com/flengyel/Dosen-Petric-Maximality/tree/dosen-petric-v25-palomar-v3},
}

@misc{Notation,
  author = {Florian Lengyel},
  title = {{The Unbearable Notation of Logic}},
  howpublished = {Version 8, Section 9},
  month = {September 19,},
  year = {2026},
  url = {https://github.com/flengyel/notes/blob/42eaccdccf6413fe4d643a347fab33edafb1c056/logic/The-Unbearable-Notation-of-Logic_v8.pdf},
}

@article{DPCCC,
  author = {Kosta Do{\v{s}}en and Zoran Petri{\'c}},
  title = {{The maximality of the typed lambda calculus and of cartesian closed categories}},
  journal = {Publications de l'Institut Math{\'e}matique (N.S.)},
  volume = {68(82)},
  year = {2000},
  pages = {1--19},
  note = {Author version: arXiv:math/9911073v2, September 25, 2012},
  url = {https://arxiv.org/abs/math/9911073},
}
\endgroup
\end{document}